%% file: bgm4.tex
\documentclass[10pt]{amsart}
\usepackage{a4,amssymb,amsfonts,verbatim,pstricks,pst-plot,xypic,times}
\usepackage{yfonts}

\newcommand{\RR}{{\mathbb R}}
\newcommand{\ZZ}{{\mathbb Z}}

\newcommand{\MM}{{\mathbb M}}

\newcommand{\cZ}{{\mathcal Z}}
\newcommand{\cG}{{\mathcal G}}
\newcommand{\cS}{{\mathcal S}}
\newcommand{\cF}{{\mathcal F}}

\newcommand{\asin}{\mathfrak s} 
\newcommand{\acos}{\mathfrak c} 
\newcommand{\n}{\nabla}
\newcommand{\m}{^{-1}}
\renewcommand{\phi}{\varphi}
\newcommand{\phid}{\dot{\varphi}}
\newcommand{\eps}{\varepsilon}
\newcommand{\half}{\frac{1}{2}}
\renewcommand{\div}{\operatorname{div}}
\newcommand{\divM}{\operatorname{div}^M}
\newcommand{\divMt}{\operatorname{div}^{M_t}}
\newcommand{\grad}{\operatorname{grad}}
\newcommand{\gradM}{\operatorname{grad}^M}

\renewcommand{\Re}{\operatorname{Re}}

\newcommand{\scal}{\operatorname{Scal}}
\newcommand{\scalZ}{\operatorname{Scal}^\cZ}
\newcommand{\scalMt}{\operatorname{Scal}^{M_t}}
\newcommand{\tr}{\operatorname{tr}}
\newcommand{\<}{\left\langle}       
\renewcommand{\>}{\right\rangle}       
\newcommand{\pa}[2]{\frac{\partial #1}{\partial #2}}
\newcommand{\SZ}{{\Sigma \cZ}}
\newcommand{\SX}{{\Sigma X}}

\newcommand{\SpX}{{\Sigma_p X}}
\newcommand{\SM}{{\Sigma M}}

\newcommand{\nZ}{\nabla^\cZ}
\newcommand{\nSZ}{\nabla^{\Sigma \cZ}}
\newcommand{\nM}{\nabla^M}
\newcommand{\nMt}{\nabla^{M_t}}
\newcommand{\nSM}{\nabla^{\Sigma M}}
\newcommand{\nX}{\nabla^X}
\newcommand{\nSX}{\nabla^{\Sigma X}}
\newcommand{\oX}{\omega^{X}}
\newcommand{\oSX}{\omega^{{\Sigma X}}}
\newcommand{\RX}{R^{X}}
\newcommand{\RSX}{R^{\Sigma X}}
\newcommand{\RZ}{R^{\cZ}}
\newcommand{\RZk}{R^{\cZ}_\kappa}
\newcommand{\RSZ}{R^{\Sigma \cZ}}
\newcommand{\RicZ}{{\operatorname{Ric}}^\cZ}
\newcommand{\ricZ}{{\operatorname{ric}}^\cZ}
\newcommand{\RicX}{{\operatorname{Ric}}^X}
\newcommand{\ricX}{{\operatorname{ric}}^X}
\newcommand{\ricMt}{{\operatorname{ric}}^{M_t}}
\newcommand{\DDM}{\tilde{D}^{M}}
\newcommand{\PSOM}{P_{\mathrm{SO}}(M)}
\newcommand{\PSOX}{P_{\mathrm{SO}}(X)}
\newcommand{\PSOZ}{P_{\mathrm{SO}}(\cZ)}

\newcommand{\PSpinM}{P_{\mathrm{Spin}}(M)}
\newcommand{\PSpinZ}{P_{\mathrm{Spin}}(\cZ)}
\newcommand{\PSpinX}{P_{\mathrm{Spin}}(X)}

\newcommand{\Ors}{{\operatorname{O}(r,s)}} 
\newcommand{\SOrs}{{\operatorname{SO}(r,s)}} 
\newcommand{\SO}{{\operatorname{SO}}} 
\newcommand{\SL}{{\operatorname{SL}}} 
\newcommand{\GL}{{\operatorname{GL}}} 
\newcommand{\GLt}{\widetilde{\GL}} 
\newcommand{\Clrs}{{\operatorname{Cl}_{r,s}}}
\newcommand{\Clors}{{\operatorname{Cl}}^0_{r,s}}
\newcommand{\Cllrs}{{\operatorname{Cl}}_{r,s}^1}
\newcommand{\Clrls}{{\operatorname{Cl}}_{r+1,s}}
\newcommand{\Clorls}{{\operatorname{Cl}}^0_{r+1,s}}
\newcommand{\Spinrs}{{\operatorname{Spin}(r,s)}} 
\newcommand{\Ad}{{\operatorname{Ad}}}
\newcommand{\Srs}{\Sigma_{r,s}}
\newcommand{\Srls}{\Sigma_{r+1,s}}
\newcommand{\Sprls}{\Sigma_{r+1,s}^+}
\newcommand{\Smrls}{\Sigma_{r+1,s}^-}
\newcommand{\Sprs}{\Sigma_{r,s}^+}
\newcommand{\Smrs}{\Sigma_{r,s}^-}
\newcommand{\Sors}{\Sigma_{r,s}^0}
\newcommand{\Slrs}{\Sigma_{r,s}^1}
\newcommand{\Aut}{{\operatorname{Aut}}} 
\newcommand{\id}{{\operatorname{id}}} 
\newcommand{\vol}{{\operatorname{vol}}} 
\newcommand{\DW}{\mathfrak{D}^W}
\newcommand{\summe}{\sum_{i=1}^n}
\newcommand{\gdt}{\dot{g}_t}
\newcommand{\gddt}{\ddot{g}_t}
\newcommand{\fd}{\dot{f}}
\newcommand{\fdd}{\ddot{f}}

\newcommand{\cdott}{\bullet_t}
\newcommand{\Mt}{{M_t}}
\newcommand{\Mto}{{M_{t_0}}}
\newcommand{\Mtl}{{M_{t_1}}}
\newcommand{\ttt}{{\tau_{t_0}^{t_1}}}
\newcommand{\DDMt}{{D}^{\Mt}}
\newcommand{\gradMt}{\operatorname{grad}^{\Mt}}
\newcommand{\DWt}{\mathfrak{D}^{W_t}}
\newcommand{\Dgdt}{\mathfrak{D}^{\gdt}}
\newcommand{\hml}{^{-1}}
\newcommand{\dichtM}{\Omega^{|n|}(M)}

\newtheorem{thm}{Theorem}[section]
\newtheorem{lemma}[thm]{Lemma}
\newtheorem{prop}[thm]{Proposition}
\newtheorem{cor}[thm]{Corollary}
\newtheorem{remark}[thm]{Remark}
\newtheorem{remarks}[thm]{Remarks}
\newtheorem{definition}[thm]{Definition}
\newtheorem{notation}[thm]{Notation}
\newtheorem{example}[thm]{Example}

\newcommand{\Example}[1]{\begin{example}{\rm #1}\end{example}}

\begin{document}
%%%%%%%%%%%%%%%%%%%%%%%%%%%%%%%%%%%%%%%%%%%%%%%%%%%%%%%%%%%%%%%%%%%%%%%%%

\title
[Generalized Cylinders in Semi-Riemannian and Spin Geometry]
{Generalized Cylinders in Semi-Riemannian and Spin Geometry}

\author{Christian B{\"a}r}
\author{Paul Gauduchon}
\author{Andrei Moroianu}

\address{
Universit{\"a}t Hamburg\\
FB Mathematik\\
Bundesstr.~55\\
20146 Hamburg\\
Germany
}
\address{
{\'E}cole Polytechnique\\
Centre de Math{\'e}matiques\\
91128 Palaiseau Cedex\\
France
}
\address{
{\'E}cole Polytechnique\\
Centre de Math{\'e}matiques\\
91128 Palaiseau Cedex\\
France
}
 
\email{baer@math.uni-hamburg.de} 
\email{pg@math.polytechnique.fr}
\email{am@math.polytechnique.fr}



\thanks{}

\date{\today}

\begin{abstract}
We use a construction which we call generalized cylinders to give a new
proof of the fundamental theorem of hypersurface theory.
It has the advantage of being very simple and the result directly extends
to semi-Riemannian manifolds and to embeddings into spaces of constant 
curvature.
We also give a new way to identify spinors for different metrics and
to derive the variation formula for the Dirac operator. 
Moreover, we show that generalized Killing spinors for Codazzi tensors
are restrictions of parallel spinors.
Finally, we study the space of Lorentzian metrics and give a criterion
when two Lorentzian metrics on a manifold can be joined in a natural
manner by a 1-parameter family of such metrics.
\end{abstract}

\maketitle

%%%%%%%%%%%%%%%%%%%%%%%%%%%%%%%%%%%%%%%%%%%%%%%%%%%%%%%%%%%%%%%%%%%%%%%%%
\section{Introduction}
%%%%%%%%%%%%%%%%%%%%%%%%%%%%%%%%%%%%%%%%%%%%%%%%%%%%%%%%%%%%%%%%%%%%%%%%%

\yinipar{I}n this paper we give various applications of a construction
which we call {\em generalized cylinders}. 
Let $M$ be a manifold and let $g_t$ be a smooth 1-parameter family of
semi-Riemannian metrics on $M$, $t\in I \subset \RR$.
Then we call the manifold $\cZ = I \times M$ with the metric $dt^2 + g_t$
a generalized cylinder over $M$.
On the one hand, this ansatz is very flexible.
Locally, near a semi-Riemannian hypersurface with spacelike normal bundle
every semi-Riemannian manifold is of this form.
The restriction to spacelike normal bundle, i.~e.\ to the positive sign
in front of $dt^2$ in the metric of $\cZ$ is made for convenience only.
Changing the signs of the metrics on $M$ as well as on $\cZ$ reduces the
case of a timelike normal bundle to that of a spacelike normal bundle.
On the other hand, this ansatz still allows to closely relate the geometries
of $M$ and $\cZ$.

In Section~\ref{secspingeo} we collect basic material on spinors and the
Dirac operator on semi-Riemannian manifolds.
We do this to fix notation and for the convenience of the reader.
Some of the material, such as the spin geometry of submanifolds, is not
so easily found in the literature unless one restricts oneself to the
Riemannian situation.

In Section~\ref{secfoliation} we study spinors on a manifold foliated
by semi-Riemannian hypersurfaces.
In particular, we derive a formula for the commutator of the leafwise
Dirac operator and the normal derivative.
This formula will be important later.

In Section~\ref{secgenzyl} we collect formulas relating the curvature
of a generalized cylinder to geometric data on $M$.

After these preliminaries we give a first application in Section~\ref{secident}.
One technical difficulty when dealing with spinors comes from the fact
that the definition of spinors depends on the metric on the manifold.
This problem does not arise when one works with tensors.
Thus if one wants to compare the Dirac operators for two different metrics,
then one first has to identify the spinor bundles in a natural manner.
This identification problem can be split into two steps.
First, construct an identification for 1-parameter families of metrics
and, secondly, given two metrics construct a natural 1-parameter family
joining them.

The second step is trivial for Riemannian metrics; just use linear interpolation.
For indefinite semi-Riemannian metrics the situation is much more complicated.
In fact, two semi-Riemannian metrics on a manifold cannot always be joined 
by a continuous path of metrics even if they have the same signature.
In Section~\ref{seclometrics} we study this problem in detail for Lorentzian metrics
and we give a criterion when two Lorentzian metrics can be joined in a natural
manner.

The first step, identifying spinors for 1-parameter families of semi-Riemannian
metrics, is carried out in Section~\ref{secident}.
The idea is very simple.
Given a 1-parameter family of metrics take the corresponding generalized 
cylinder and use parallel transport on this cylinder.
%Some technicalities have to be worked out, such as the identification
%of spinors on the manifold and on the generalized over it.
It turns out that this identification is the same as the one constructed
differently by Bourguignon and the second author in 
\cite{bourguignon-gauduchon92a} for Riemannian metrics.
The commutator formula from Section~\ref{secfoliation} directly translates
to the variation formula for Dirac operators.

This variation formula is what one needs to compute the energy-momentum
tensor for spinors.
To make this precise we briefly summarize Lagrangian field theory in 
Section~\ref{secenergymomentum} and we give a general definition of 
energy-momentum tensors.
Then we compute the example of the Lagrangian for spinors given by the
Dirac operator.

In Section~\ref{sechyper} we give a new and simple proof of the fundamental theorem
of hypersurface theory.
A hypersurface of $\RR^{n+1}$ inherits a Riemannian metric and its Weingarten
map must satify the Gauss and Codazzi-Mainardi equations.
The fundamental theorem says that, conversely, any Riemannian manifold $M$
with a symmetric endomorphism field of $TM$ satisfying the Gauss and
Codazzi-Mainardi equations can, at least locally, be embedded isometrically
into $\RR^{n+1}$ with Weingarten map given by this endomorphism field.
Our proof goes like this:
We write down an {\em explicit} metric on the cylinder $\cZ = I \times M$
and we then check that this metric is flat.
Since every flat Riemannian manifold is locally isometric to Euclidean space
the theorem follows.
This approach directly extends to semi-Riemannian manifolds and to embeddings
into spaces of constant sectional curvature not necessarily zero.
This kind of approach to the fundamental theorem for hypersurfaces
was suggested, but not carried out, by Petersen in \cite[p.~95]{petersen98a}.

In Section~\ref{seckilling} we study generalized Killing spinors.
They are characterized by the overdetermined equation 
$\n^\SM_X\psi=\frac{1}{2}A(X)\cdot\psi$ where $A$ is a given symmetric
endomorphism field.
We show that if $A$ is a Codazzi tensor, then the manifold can be embedded
as a hypersurface into a Ricci flat manifold equipped with a parallel
spinor which restricts to $\psi$.
This generalizes the case of Killing spinors, $A = \lambda\,\id$.
The classification of manifolds admitting Killing spinors in \cite{baer93a}
was based on the observation that the cone over such a manifold possesses
a parallel spinor.
This also generalizes the case that $A$ is parallel which was studied in
\cite{morel03a}.

{\bf Acknowledgements.}
The authors would like to thank W.~Ballmann and H.~Karcher for valuable
suggestions.
The authors have been partially supported by the Research and Training 
Network HPRN-CT-2000-00101 ``EDGE'' funded by the European Commission.
The first author has also been partially supported by the Research and Training 
Network HPRN-CT-1999-00118 ``Geometric Analysis''.
The first author would like to thank the Ecole Polytechnique, Palaiseau, and the 
the Max-Planck-Institut f\"ur Mathematik, Bonn, for their hospitality.

%%%%%%%%%%%%%%%%%%%%%%%%%%%%%%%%%%%%%%%%%%%%%%%%%%%%%%%%%%%%%%%%%%%%%%%%%
\section{The Dirac operator on semi-Riemannian manifolds}
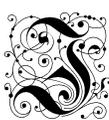
\label{secspingeo}
%%%%%%%%%%%%%%%%%%%%%%%%%%%%%%%%%%%%%%%%%%%%%%%%%%%%%%%%%%%%%%%%%%%%%%%%%

\yinipar{I}n this section we collect the basic facts and conventions concerning
spinors and Dirac operators on semi-Riemannian manifolds.
For a detailed introduction the reader may consult the book \cite{baum81a}.
We start with some algebraic preliminaries.
Let $r+s=n$ and consider the nondegenerate symmetric bilinear form of signature
$(r,s)$
$$
\< v,w \> := \sum_{i=1}^r v^i w^i - \sum_{i=r+1}^n v^i w^i
$$
on $\RR^n$.
Define the corresponding {\em orthogonal group} by
$$
\Ors := \{ A\in {\rm GL}(n,\RR)\ |\ \< Av,Aw \> = \< v,w \> 
\mbox{ for all } v,w \in \RR^n  \}
$$
and the {\em special orthogonal group} by
$$
\SOrs := \{ A\in \Ors\ |\ \det(A) = 1 \}.
$$
If $r=0$ or $s=0$, then $\SOrs$ is connected, otherwise it has 
two connected components.
%This is due to the fact that for $r,s \ge 1$ there are both
%timelike vectors satisfying $\<v,v\> < 0$ as well as nonzero spacelike vectors
%satisfying $\<v,v\> > 0$.
%Both these sets are acted upon by $\Ors$ and they carry natural orientations.
%The connected component of the identity of the group $\SOrs$ (or,
%equivalently, of the group $\Ors$) is the set of those $A \in \Ors$ which
%respect both these orientations.
%The elements of the second component of $\SOrs$ invert both orientations.

Now let $\Clrs$ be the {\em Clifford algebra} corresponding to the symmetric 
bilinear form $\<\cdot,\cdot\>$.
This is the unital algebra generated by $\RR^n$ subject to the relations
\begin{equation}
v\cdot w + w\cdot v + 2 \<v,w\> \cdot 1 = 0
\label{cliffrel}
\end{equation}
for all $v,w\in \RR^n$.
There is a decomposition into even and odd elements
$$
\Clrs = \Clors \oplus \Cllrs
$$ 
such that $\RR$ injects naturally into $\Clors$ and $\RR^n$ into 
$\Cllrs$.
The {\em spin group} is defined by
$$
\Spinrs := \{ v_1 \cdots v_k \in \Clors\ |\ v_j \in \RR^n 
\mbox{ such that }\<v_j,v_j\> = \pm 1 
\mbox{ and $k$ is even}\}
$$
with multiplication inherited from $\Clrs$.
Given $v\in\RR^n$ such that $\<v,v\>\not=0$ and arbitrary $w\in\RR^n$
we see directly from relation (\ref{cliffrel}) that $v^{-1} =
-\frac{v}{\<v,v\>}$ and
$$
\Ad_v(w) := v^{-1}\cdot w\cdot v = - w + 2 \frac{\<v,w\>}{\<v,v\>}v.
$$
Hence $-\Ad_v$ is the reflection across the hyperplane $v^\perp$ and,
in particular, leaves $\RR^n \subset \Clrs$ invariant.
Thus conjugation gives an action of $\Spinrs$ on $\RR^n$ by an even
number of reflections across hyperplanes.
This yields the exact sequence
$$
1 \longrightarrow \ZZ/2\ZZ = \{1,-1\}  \longrightarrow \Spinrs
\stackrel{\Ad}{\longrightarrow}  \SOrs  \longrightarrow 1.
$$
If $n=r+s$ is even the Clifford algebra possesses an irreducible complex 
module $\Sigma_{r,s}$ of complex dimension dimension $2^{n/2}$, the complex
{\em spinor module}.
When restricted to $\Clors$ the spinor module decomposes into
$$
\Srs = \Sprs \oplus \Smrs ,
$$
the submodules of spinors of {\em positive} resp.\ {\em negative chirality}.
In particular, the spin group $\Spinrs \subset \Clors$ acts on $\Sprs$
and on $\Smrs$.
This action
$$
\rho = \rho^+ \oplus \rho^- : \Spinrs \to \Aut(\Sprs) \times \Aut(\Smrs)
\subset \Aut(\Srs)
$$
is called the {\em spinor representation} of $\Spinrs$.
Given an orientation on $\RR^n$ the $\Clors$-submodules $\Sprs$ and $\Smrs$ can be 
characterized by the action of the volume element $\vol := e_1 \cdots e_n \in \Clors$
which acts on $\Sprs$ as $+i^{s+n(n+1)/2}\id$ and on $\Smrs$ as 
$-i^{s+n(n+1)/2}\id$ where $e_1,\ldots,e_n$ is a positively oriented orthonormal 
basis of $\RR^n$.

If $n$ is odd, then $\Clrs$ has two inequivalent irreducible modules
$\Sors$ and $\Slrs$, both of complex dimension $2^{(n-1)/2}$.
These two modules are again distinguished by the action of the volume
element $\vol = e_1 \cdots e_n \in \Cllrs$, namely $\vol$ acts
as $+i^{s+n(n+1)/2}\id$ on $\Sors$ and as $-i^{s+n(n+1)/2}\id$ on $\Slrs$.
When restricted to $\Clors$ the two modules become equivalent and we
simply write $\Srs := \Sors$.
This time the spinor representation
$$
\rho: \Spinrs \to \Aut(\Srs)
$$
is irreducible.
All spinor modules carry nondegenerate symmetric sesquilinear forms $\<\cdot,
\cdot\>$ (in general not definite) which are invariant under the
action of $\Spinrs$.
The action of a vector $v \in \RR^n \subset \Clrs$ on $\Srs$ is skewsymmetric
with respect to $\<\cdot,\cdot\>$, i.~e.\ $\<v\cdot\sigma_1,\sigma_2\>
= -\<\sigma_1,v\cdot\sigma_2\>$.

To prepare for the study of submanifolds later on we now look at an
embedding of $\RR^n$ into $\RR^{n+1}$ such that $(\RR^n)^\perp$ is
spacelike.
Let $(\RR^n)^\perp$ be spanned by a spacelike unit vector $e_0$.
The map $\RR^n \to \Clrls$, $v \mapsto e_0\cdot v$, induces an algebra
isomorphism $\Clrs \to \Clorls$ under which the volume element of 
$\Clrs$ is mapped to the volume element of $\Clrls$ in case $n$ is odd.

If $n$ is even, then $\Srls$ pulls back to $\Srs$ under this algebra
isomorphism.
In other words, we can regard $\Srls$ as the spinor representation of
$\Clrs$ provided we define the action of $\Clrs$ on $\Srls$ by
$$
v \otimes \sigma \mapsto e_0 \cdot v \cdot \sigma
$$
where $v\in \RR^n$ and $\cdot$ denotes the action of $\Clrls$.

Similarly, if $n$ is odd, then the action of the volume forms
shows that $\Sprls$ pulls back to $\Sors$ while 
$\Smrls$ pulls back to $\Slrs$.

Now we turn to geometry.
Let $X$ denote an oriented $n$-dimensional differentiable manifold.
The bundle $P_{\GL^+}(X)$ of positively oriented tangent frames
forms a $\GL^+(n,\RR)$-principal bundle over $X$.
Here and henceforth $\GL^+(n,\RR)$ denotes the group of real 
$n\times n$-matrices with positive determinante and $A:\GLt^+(n,\RR)\to 
\GL^+(n,\RR)$ its connected twofold covering group.
A {\em spin structure} of $X$ is a $\GLt^+(n,\RR)$-principal bundle
$P_{\GLt}(X)$ over $X$ together with a twofold covering map 
$\Theta:P_{\GLt^+}(X) \to P_{\GL^+}(X)$ such that the following diagram 
commutes
\begin{equation}
\xymatrix{
P_{\GLt^+}(X)\times \GLt^+(n,\RR) \ar[dd]^{\Theta \times A} \ar[r] & P_{\GLt^+}(X) \ar[dd]_\Theta \ar[dr] & \\
& & X\\
P_{\GL^+}(X)\times \GL^+(n,\RR)  \ar[r]   & P_{\GL^+}(X) \ar[ur] & 
}
\label{spindiagram}
\end{equation}
where the horizontal arrows denote the group actions on the principal bundles.
This definition of a spin structure has the advantage of being independent 
of the choice of any semi-Riemannian metric on $X$.
An oriented manifold together with a spin structure will be called a 
{\em spin manifold}.

Let $X$ now in addition carry a semi-Riemannian metric of signature $(r,s)$, 
$r+s=n$.
The bundle $\PSOX\subset P_{\GL^+}(X)$ of positively oriented {\em orthonormal}
tangent frames forms an $\SOrs$-principal bundle over $X$.
Restricting $A:\GLt^+(n,\RR)\to \GL^+(n,\RR)$ to the preimage of $\SOrs
\subset\GL^+(n,\RR)$ we recover $\Ad : \Spinrs \to \SOrs$.
Putting $\PSpinX := \Theta\hml(\PSOX)$ we get a $\Spinrs$-principal bundle and
and the maps in diagram (\ref{spindiagram}) restrict to the following 
commutative diagram
$$
\xymatrix{
\PSpinX \times \Spinrs \ar[dd]^{\Theta \times \Ad} \ar[r] & \PSpinX \ar[dd]_\Theta \ar[dr] & \\
& & X\\
\PSOX \times \SOrs  \ar[r]   & \PSOX \ar[ur] & 
}
$$
Very often in the literature $\PSpinX$ is called a spin structure of $X$ and 
we will call $X$ together with $\PSpinX$ a {\em semi-Riemannian spin manifold}.

On a semi-Riemannian spin manifold we define the {\em spinor bundle}
of $X$ as the complex vector bundle associated to the spinor representation,
i.~e.\ 
$$
\SX := \PSpinX \times_\rho \Srs.
$$
In other words, for $p\in X$ the fiber of $\SpX$ of $\SX$ over $p$ consists
of equivalence classes of pairs $[b,\sigma]$ where $b\in \PSpinX_p$ and
$\sigma\in\Srs$ subject to the relation
$$
[b,\sigma] = [bg^{-1},g\sigma]
$$
for all $g\in \Spinrs$. 
Unfortunately, the spinor bundle cannot be defined independently of the metric
using $P_{\GLt^+}(X)$ instead of $\PSpinX$ because the spinor representation
$\rho$ of $\Spinrs$ on $\Srs$ does not extend to a representation of 
$\GLt^+(n,\RR)$ on $\Srs$.
We will come back to this problem in Section~\ref{secident}.

Note that the tangent bundle can also be written in a similar manner,
$TX = \PSOX \times_\tau \RR^n$ where $\tau$ is the standard representation
of $\SO(r,s)$ on $\RR^n$.
One defines {\em Clifford multiplication} $T_pX \otimes \SpX \to \SpX$
by 
$$
[\Theta(b),v] \cdot [b,\sigma] := [b,v\cdot\sigma] 
$$
where $b\in \PSpinX_p$, $v\in\RR^n$, and $\sigma\in\Srs$.
For $g\in \Spinrs$ we see from
\begin{eqnarray*}
[\Theta(bg),v]\cdot[bg,\sigma] &=& [\Theta(b)\Ad_g,v]\cdot[bg,\sigma]
\, = \,  [\Theta(b),\Ad_g v]\cdot[b,g\sigma] \\
&=& [b,gvg^{-1}g\sigma] \, = \,  [b,gv\sigma] \, = \,  [bg,v\sigma]
\end{eqnarray*}
that this is well-defined.
It is this point that goes wrong when one tries to work with 
nonoriented manifolds and pin structures.
Had we defined $\Srs = \Slrs$ instead of $\Srs = \Sors$ in odd dimensions,
then we would have obtained the Clifford multiplication with the opposite
sign.

Clifford multiplication inherits the relations of the Clifford algebra,
i.~e.\ for $X,Y \in T_pX$ and $\phi\in\SpX$ we have
$$
X\cdot Y \cdot \phi + Y \cdot X \cdot \phi + 2 \<X,Y\> \phi = 0.
$$

In even dimensions the spinor bundle splits into the positive and the
negative {\em half-spinor bundles}, 
\begin{equation}
\SX = \Sigma^+X \oplus \Sigma^-X
\label{spinsplit}
\end{equation}
where $\Sigma^\pm X = \PSpinX \times_{\rho^\pm} \Sigma^\pm_{r,s}$.
Clifford multiplication by a tangent vector interchanges $\Sigma^+X$
and $\Sigma^-X$.

The $\Spinrs$-invariant nondegenerate symmetric sesquilinear forms on $\Srs$
and $\Sigma^\pm_{r,s}$ induce (in general indefinite) inner products
on $\SX$ and $\Sigma^\pm X$ which we again denote by $\<\cdot,\cdot\>$.

The connection 1-form $\oX$ on $\PSOX$ for the Levi-Civita connection $\nX$ 
can be lifted via $\Theta$ to $\PSpinX$, i.~e.\ $\oSX := 
\Ad_*^{-1}\circ\Theta^*(\oX)$.
Composing with $\Ad_*^{-1}$ is necessary because the connection 1-form 
on $\PSpinX$ must take values in the Lie algebra of $\Spinrs$ rather
than in that of $\SOrs$.
Now $\oSX$ induces a covariant derivative $\nSX$ on $\SX$.

An equivalent, but less invariant, way of describing $\nSX$ is as follows:
If $b$ is a local section in $\PSpinX$, then $\Theta(b) = (e_1,\ldots,e_n)$
is a local oriented orthonormal tangent frame, $\<e_i,e_j\> \equiv
\eps_i \delta_{ij}$ where $\eps_i = \pm 1$.
The Christoffel symbols of $\nX$ with respect to this frame are given by
$$
\nX_{e_i}e_j = \sum_{k=1}^n \Gamma_{ij}^k\, e_k .
$$
Now the covariant derivative of a locally defined spinor field
$\phi = [b,\sigma]$, $\sigma$ a function with values in $\Srs$, is given by
\begin{equation}
\nSX_{e_i}\phi = 
\left[b, d_{e_i}\sigma + 
\frac{1}{2}\sum_{j<k}\Gamma_{ij}^k\, \eps_j\, e_j \cdot e_k \cdot \sigma\right] .
\label{zshg}
\end{equation}

One checks that $\nSX$ is a metric connection and that it leaves 
the splitting (\ref{spinsplit}) in even dimensions invariant.
Moreover, it satisfies the following Leibniz rule:
$$
\nSX_Z(Y\cdot \phi) = (\nX_ZY)\cdot\phi + Y\cdot\nSX_Z \phi
$$
for all vector fields $Z$ and $Y$ and all spinor fields $\phi$.

The curvature tensor $\RSX$ of $\nSX$ can be computed in terms of
the curvature tensor $\RX$ of the Levi-Civita connection,
$$
\RSX(Y,Z)\phi = 
\frac12 \sum_{i<j} \eps_i\eps_j\< \RX(Y,Z)e_i,e_j\>e_i \cdot e_j \cdot\phi .
$$
Using the first Bianchi identity one easily computes
\begin{equation}
\sum_{i=1}^n \eps_i\, e_i \cdot \RSX(e_i,Y)\phi =
\half \RicX(Y) \cdot \phi .
\label{ricci}
\end{equation}
Here $\RicX$ denotes the {\em Ricci curvature} considered as an endomorphism
field on $TM$.
The Ricci curvature considered as a symmetric bilinear form will be
written $\ricX(Y,Z) = \<\RicX(Y),Z\>$.

The {\em Dirac operator} maps spinor fields to spinor fields and is defined
by
$$
D^X\phi = \sum_{i=1}^n \eps_i e_i \cdot \nSX_{e_i} \phi .
$$
Given two spinor fields $\phi$ and $\psi$ one can define a vector
field $Y$ by the requirement $\<Y,Z\> = \<Z\cdot\phi,\psi\>$ for all vector
fields $Z$ and one easily computes
$$
\div(Y) = \<D^X\phi,\psi\> - \<\phi,D^X\psi\>.
$$
Hence the Dirac operator is formally selfadjoint, i.~e.\ if the 
intersection of the supports of $\phi$ and $\psi$ is compact, then
$$
(D^X\phi,\psi) = (\phi,D^X\psi)
$$
where $(\phi,\psi) = \int_M \<\phi,\psi\> dV$.

%%%%%%%%%%%%%%%%%%%%%%%%%%%%%%%%%%%%%%%%%%%%%%%%%%%%%%%%%%%%%%%%%%%%%%%%%
\section{The Dirac operator on manifolds foliated by hypersurfaces}
\label{secfoliation} 
%%%%%%%%%%%%%%%%%%%%%%%%%%%%%%%%%%%%%%%%%%%%%%%%%%%%%%%%%%%%%%%%%%%%%%%%%

\yinipar{L}et $\cZ$ be an oriented $(n+1)$-dimensional semi-Riemannian spin manifold.
Let $\Theta : \PSpinZ \to \PSOZ$ be a spin structure on $\cZ$.
Let $M \subset \cZ$ be a semi-Riemannian hypersurface with trivial 
spacelike normal bundle.
This means there is a vector field $\nu$ on $\cZ$ along $M$ satisfying
$\<\nu,\nu\> = +1$ and $\<\nu,TM\>=0$.
If the signature of $M$ is $(r,s)$, then the signature of $\cZ$
is $(r+1,s)$.

In this situation $M$ inherits a spin structure as follows:
The bundle of oriented orthonormal frames of $M$, $\PSOM$, can be embedded
into the bundle of oriented orthonormal frames of $\cZ$ restricted to $M$, 
$\PSOZ|_M$, by the map $\iota: (e_1,\ldots,e_n) \mapsto (\nu,e_1,\ldots,e_n)$.
Then $\PSpinM := \Theta\m(\iota(\PSOM))$ defines a spin structure on $M$.
We will always implicitly assume that this spin structure be taken on $M$.
The same discussion is possible on the level of $\GLt^+(n,\RR)$-bundles.

The algebraic remarks in the previous section show that if
$n$ is even, then 
$$
\SZ|_M = \SM
$$
where the Clifford multiplication with respect to $M$ is given by
$X \otimes \phi \mapsto \nu \cdot X \cdot \phi$ and ``$\cdot$''
always denotes the Clifford multiplication with respect to $\cZ$.
If $n$ is odd, then
$$
\Sigma^+\cZ|_M = \SM
$$
and again Clifford multiplication with respect to $M$ is given
by $X \otimes \phi \mapsto \nu \cdot X \cdot \phi$ while
$$
\Sigma^-\cZ|_M = \SM
$$
with Clifford multiplication with respect to $M$ given
by $X \otimes \phi \mapsto  - \nu \cdot X \cdot \phi$.
The minus sign comes from the fact that in odd dimensions we defined
$\Srs = \Sors$ while $\Slrs$ leads to the opposite sign for the
Clifford multiplication.
The identifications preserve the natural inner products $\<\cdot,\cdot\>$.

Let $W$ denote the {\em Weingarten map} with respect to $\nu$, i.~e.\
\begin{equation}
\nZ_XY = \nM_XY + \<W(X),Y\>\nu
\label{weingarten}
\end{equation}
for all vector fields $X$ and $Y$ on $M$.
The Weingarten map is symmetric with respect to the semi-Riemannian
metric, $\<W(X),Y\> = \<X,W(Y)\>$ and is also given by $W(X) = -\nZ_X\nu$.
If we denote the Christoffel symbols of $M$ with respect to a local
orthogonal tangent frame $(e_1,\ldots,e_n)$ by $\Gamma^{M,k}_{ij}$
and the Christoffel symbols of $\cZ$ with respect to $(e_0,e_1,\ldots,e_n)$,
$e_0=\nu$, by $\Gamma^{\cZ,k}_{ij}$, then (\ref{weingarten}) implies for 
$1 \leq i,j,k \leq n$
\begin{eqnarray}
\Gamma^{\cZ,k}_{ij} &=& \Gamma^{M,k}_{ij} ,\\
\Gamma^{\cZ,0}_{ij} &=& \<W(e_i),e_j\> ,\\
\Gamma^{\cZ,k}_{i0} &=& -\eps_0\eps_k \Gamma^{\cZ,0}_{ik} 
                   = -\eps_k\<W(e_i),e_k\> 
\label{Gammazyl0}.
\end{eqnarray}
Plugging this into (\ref{zshg}) we get for a section $\phi = [b,\sigma]$
of $\SZ|_M$ and $1 \leq i \leq n$
\begin{eqnarray*}
\nSZ_{e_i}\phi &=&
\left[b, d_{e_i}\sigma + 
\half\left(-\sum_{k=1}^n\eps_k\<W(e_i),e_k\>\eps_0 e_0\cdot e_k + 
\sum_{1\leq j < k\leq n}\Gamma^{M,k}_{ij}\eps_je_j\cdot e_k\right)
\cdot \sigma\right] \\
&=&
\left[b, d_{e_i}\sigma + 
\half\left(- e_0\cdot W(e_i)  + 
\sum_{1\leq j < k\leq n}\Gamma^{M,k}_{ij}\eps_je_0\cdot e_j\cdot e_0\cdot e_k\right)
\cdot \sigma\right] \\
&=&
\nSM_{e_i}\phi - \half \nu \cdot W(e_i) \cdot \phi .
\end{eqnarray*}
Hence for each $X\in TM$ and each section $\phi$ of $\SZ|_M$ we have
\begin{equation}
\nSZ_{X}\phi  =  \nSM_{X}\phi - \half \nu \cdot W(X) \cdot \phi .
\label{spingauss}
\end{equation}
Now let $\phi$ be a section of $\SZ$ defined in a neighborhood of $M$.
On the one hand,
$$
D^\cZ\phi = \sum_{i=1}^n \eps_i e_i \cdot \nSZ_{e_i} \phi 
+ \nu \cdot \nSZ_{\nu} \phi .
$$
On the other hand by (\ref{spingauss}),
\begin{eqnarray*}
\sum_{i=1}^n \eps_i e_i \cdot \nSZ_{e_i} \phi 
&=&
\sum_{i=1}^n \eps_i\, e_i \cdot \nSM_{e_i} \phi 
-\half \sum_{i=1}^n \eps_i\, e_i \cdot \nu \cdot W(e_i) \cdot \phi \\
&=&
- \nu\cdot \sum_{i=1}^n \eps_i\, \nu\cdot e_i \cdot \nSM_{e_i} \phi 
+\half \sum_{i=1}^n \eps_i\, \nu \cdot e_i \cdot  W(e_i) \cdot \phi \\
&=&
- \nu\cdot \DDM - \half \tr(W) \nu \cdot\phi
\end{eqnarray*}
where $\DDM = D^M$ if $n$ is even and $\DDM = 
\begin{pmatrix}
D^M & 0 \cr 
0 & -D^M
\end{pmatrix}$
if $n$ is odd.
Thus the Dirac operators on $M$ and on $\cZ$ are related by
\begin{equation}
\nu\cdot D^\cZ = \DDM +\frac{n}{2}H - \nSZ_\nu
\label{diracgauss}
\end{equation}
where $H = \frac1n \tr(W)$ denotes the mean curvature.

Next we consider the situation that $\cZ$ carries a semi-Riemannian foliation by 
hypersurfaces.
The commutator of the leafwise Dirac operator and the normal derivative
will be of central importance later.

\begin{prop}
Let $\cZ$ be an $(n+1)$-dimensional semi-Riemannian spin manifold.
Let $\cZ$ carry a semi-Riemannian foliation by hypersurfaces with trivial
spacelike normal bundle, i.~e.\ the leaves $M$ are semi-Riemannian 
hypersurfaces and there exists a vector field $\nu$ on $\cZ$ perpendicular
to the leaves such that $\<\nu,\nu\>=1$ and $\nZ_\nu\nu = 0$.
Let $W$ denote the Weingarten map of the leaves with respect to $\nu$
and let $H=\frac1n\tr(W)$ be the mean curvature.

Then the commutator of the leafwise Dirac operator and the normal
derivative is given by
$$
[\nSZ_\nu,\DDM]\,\phi = 
\DW\phi - \frac{n}{2}\, \nu\cdot\gradM(H)\cdot\phi 
+ \half\,\nu\cdot\divM(W)\cdot\phi  .
$$
Here $\gradM$ denotes the leafwise gradient, $\divM(W) = \sum_{i=1}^n
\eps_i\, (\nM_{e_i}W)(e_i)$ denotes the leafwise divergence of the
endomorphism field $W$, $\DW\phi =
\sum_{i=1}^n\eps_i\, \nu\cdot e_i\cdot\nSM_{W(e_i)}\phi$, and ``$\cdot$''
denotes Clifford multiplication on $\cZ$.
\label{kommutator}
\end{prop}

\begin{proof}
We choose a local oriented orthonormal tangent frame $(e_1,\ldots,e_n)$ 
for the leaves and we may assume for simplicity that $\nZ_\nu e_i = 0$.
We compute
\begin{eqnarray}
[\nSZ_\nu,\DDM]\,\phi &=&
\summe \eps_i \left( \nSZ_{\nu}(\nu \cdot e_i\cdot \nSM_{e_i}\phi) 
-  \nu \cdot e_i\cdot \nSM_{e_i} \nSZ_{\nu}\phi  \right)
\nonumber\\
&=&
\summe \eps_i\,\nu \cdot e_i\cdot  \left( \nSZ_{\nu}\nSM_{e_i}\phi 
-  \nSM_{e_i} \nSZ_{\nu}\phi  \right)
\nonumber\\
&\stackrel{(\ref{spingauss})}{=}&
\summe \eps_i\,\nu \cdot e_i\cdot  \Big( \nSZ_{\nu}(\nSZ_{e_i}+
\half \nu\cdot W(e_i))
\nonumber\\
&&
-  (\nSZ_{e_i}+\half \nu\cdot W(e_i)) \nSZ_{\nu} \Big)\phi 
\nonumber\\
&=&
\summe \eps_i\,\nu \cdot e_i\cdot \Big(\RSZ(\nu,e_i)+\nSZ_{[\nu,e_i]}
+ \half \nu\cdot (\nZ_\nu W)(e_i)   \Big)\phi
\nonumber\\
&\stackrel{(\ref{ricci})}{=}&
-\half\nu\cdot\RicZ(\nu)\cdot\phi +
\summe \eps_i\,\nu \cdot e_i\cdot \Big(\nSZ_{W(e_i)}
+ \half \nu\cdot (\nZ_\nu W)(e_i)   \Big)\phi
\nonumber\\
&\stackrel{(\ref{spingauss})}{=}&
-\half\nu\cdot\RicZ(\nu)\cdot\phi 
\nonumber\\
&&
+ \summe \eps_i\,\nu \cdot e_i\cdot \Big(\nSM_{W(e_i)} -\half\nu\cdot W^2(e_i)
+ \half \nu\cdot (\nZ_\nu W)(e_i)   \Big)\phi
\nonumber\\
&=&
-\half\nu\cdot\RicZ(\nu)\cdot\phi +\DW\phi
\nonumber\\
&&
+ \half \summe \eps_i\, e_i\cdot \Big(-W^2(e_i)
+ (\nZ_\nu W)(e_i)   \Big)\phi .
\end{eqnarray}
The Riccati equation for the Weingarten map $(\nZ_\nu W)(X) = 
\RZ(X,\nu)\nu + W^2(X)$ yields
\begin{eqnarray}
[\nSZ_\nu,\DDM]\,\phi &=&
-\half\nu\cdot\RicZ(\nu)\cdot\phi +\DW\phi
+ \half \summe \eps_i\, e_i\cdot (\RZ(e_i,\nu)\nu)\cdot\phi
\nonumber\\
&=&
-\half\nu\cdot\RicZ(\nu)\cdot\phi +\DW\phi
+ \half \ricZ(\nu,\nu)\phi
\nonumber\\
&=&
\DW\phi -\half \summe \eps_i\, \ricZ(\nu,e_i)\, \nu\cdot e_i\cdot\phi.
\label{ausdruckmitricci}
\end{eqnarray}
The Codazzi-Mainardi equation \cite[p.~115]{oneill83a} gives for $X,Y,V
\in T_pM$
$$
\<\RZ(X,Y)V,\nu\> = \<(\nM_X W)(Y),V\> - \<(\nM_Y W)(X),V\> .
$$
Thus 
\begin{eqnarray*}
\ricZ(\nu,X) &=&
\summe \eps_i \<\RZ(X,e_i)e_i,\nu\>
\nonumber\\
&=&
\summe \eps_i \left( \<(\nM_X W)(e_i),e_i\> - \<(\nM_{e_i} W)(X),e_i\>\right)
\nonumber\\
&=&
\tr(\nM_X W) - \<\divM(W),X\> .
\end{eqnarray*}
Plugging this into (\ref{ausdruckmitricci}) we get
\begin{eqnarray*}
[\nSZ_\nu,\DDM]\,\phi 
&=&
\DW\phi -\half \summe \eps_i\, \left(\tr(\nM_{e_i} W) 
- \<\divM(W),e_i\>\right)\, \nu\cdot e_i\cdot\phi\\
&=&
\DW\phi -\half \summe \eps_i\, d_{e_i} \tr(W) \nu\cdot e_i\cdot\phi
+\half \nu\cdot\divM(W) \cdot\phi\\
&=&
\DW\phi -\frac{n}{2}\nu\cdot \gradM(H)\cdot\phi 
+\half \nu\cdot\divM(W) \cdot\phi .
\end{eqnarray*}

\end{proof}

%%%%%%%%%%%%%%%%%%%%%%%%%%%%%%%%%%%%%%%%%%%%%%%%%%%%%%%%%%%%%%%%%%%%%%%%%
\section{The generalized cylinder}
\label{secgenzyl}
%%%%%%%%%%%%%%%%%%%%%%%%%%%%%%%%%%%%%%%%%%%%%%%%%%%%%%%%%%%%%%%%%%%%%%%%%

\yinipar{L}et $M$ be an $n$-dimensional differentiable manifold, let $g_t$ be
a smooth 1-parameter family of semi-Riemannian metrics on $M$,
$t\in I$ where $I \subset \RR$ is an interval.
We define the {\em generalized cylinder} by
$$
\cZ := I \times M
$$
with semi-Riemannian metric
$$
g_\cZ := dt^2 + g_t .
$$
The generalized cylinder is an $(n+1)$-dimensional semi-Riemannian
manifold (with boundary if $I$ has boundary) of signature $(r+1,s)$ 
if the signature of $g_t$ is $(r,s)$.
The vector field $\nu := \pa{}{t}$ is spacelike of unit length and
orthogonal to the hypersurfaces $M_t := \{t\}\times M$.
Let $W$ denote the Weingarten map of $M_t$ with respect to $\nu$
and let $H$ be the mean curvature.

If $X$ is a local coordinate field on $M$, then $\<X,\nu\> = 0$ and
$[X,\nu] = 0$.
Thus
\begin{eqnarray*}
0 &=& 
d_\nu \<X,\nu\> = \<\nZ_\nu X,\nu\> + \<X,\nZ_\nu\nu\>
= 
\<\nZ_X \nu,\nu\> + \<X,\nZ_\nu\nu\>\\ 
&=&
-\<W(X),\nu\> + \<X,\nZ_\nu\nu\> = \<X,\nZ_\nu\nu\>
\end{eqnarray*}
and differentiating $\<\nu,\nu\>=1$ yields $\<\nu,\nZ_\nu\nu\>=0$.
Hence
$$
\nZ_\nu\nu = 0,
$$
i.~e.\ for $p\in M$ the curves $t\mapsto (t,p)$ are geodesics parametrized
by arclength.
So the assumptions of Proposition~\ref{kommutator} are satisfied for the
foliation $(M_t)_{t\in I}$.

Now fix $p\in M$ and $X,Y \in T_pM$.
We define the first and second derivative of $g_t$ by 
\begin{eqnarray*}
\gdt(X,Y) &:=& \frac{d}{dt}(g_t(X,Y)) ,\\
\gddt(X,Y) &:=& \frac{d^2}{dt^2}(g_t(X,Y)) .
\end{eqnarray*}

Then $\gdt$ and $\gddt$ are smooth 1-parameter families of symmetric 
$(2,0)$-tensors on $M$.

\begin{prop}
On a generalized cylinder $\cZ = I \times M$ with semi-Riemannian metric $g^\cZ = 
\<\cdot,\cdot\> = dt^2 + g_t$ the following formulas hold:
\begin{eqnarray}
\< W(X),Y\> &=& -\half \gdt(X,Y) ,
\label{zylweingarten}\\
\<R^\cZ(U,V)X,Y\> &=& \<R^{M_t}(U,V)X,Y\>
\label{zylgauss}\\ 
&& + \frac{1}{4} \left(
\gdt(U,X) \gdt(V,Y) - \gdt(U,Y) \gdt(V,X)\right) ,
\nonumber \\
\<R^\cZ(X,Y)U,\nu\> &=& \half \left((\nMt_Y\gdt)(X,U) - (\nMt_X\gdt)(Y,U)\right),
\label{zylcodazzi}\\
\<R^\cZ(X,\nu)\nu,Y\> &=& -\half\left( \gddt(X,Y) + \gdt(W(X),Y)  \right),
\label{zylriccati} \\
\ricZ(\nu,\nu) &=& \tr(W^2) -\half \tr_{g_t}(\gddt) ,
\label{zylricnn}\\
\ricZ(X,\nu) &=& d_X\tr(W) - \<\divM(W),X\> ,
\label{zylrictn}\\
\ricZ(X,Y) &=& 
\ricMt(X,Y) + 2 \<W(X),W(Y)\> 
\label{zylrictt}\\
&&
-\tr(W)\<W(X),Y\>-\half \gddt(X,Y) ,
\nonumber\\
\scalZ &=& \scalMt +\, 3 \tr(W^2) - \tr(W)^2 - \tr_{g_t}(\gddt) ,
\label{zylscal}
\end{eqnarray}
where $X,Y,U,V \in T_pM$, $p\in M$.
\label{zylformeln}
\end{prop}

\begin{proof}
To show (\ref{zylweingarten})
we extend $X$ and $Y$ to local coordinate fields on $M$ so that all
Lie brackets vanish.
Then the Koszul formula \cite[p.~61]{oneill83a} for the Levi-Civita 
connection of $\cZ$ yields
\begin{eqnarray*}
\< W(X),Y\> &=& - \< \nZ_X\nu,Y \> = 
-\half \left( d_X\<\nu,Y\> + d_\nu\<Y,X\> -
d_Y\<X,\nu\>   \right) \\
&=&
-\half d_\nu\<Y,X\>
=
-\half \pa{}{t} g_t(X,Y)
= -\half \gdt(X,Y) .
\end{eqnarray*}

Equation (\ref{zylgauss}) follows directly from (\ref{zylweingarten}) and the 
Gauss equation \cite[p.~100]{oneill83a}
\begin{eqnarray*}
\<R^\cZ(U,V)X,Y\> &=& \<R^{M_t}(U,V)X,Y\> 
+ \<W(U),X\>\<W(V),Y\>\\&& - \<W(U),Y\>\<W(V),X\> .
\end{eqnarray*}

Equation (\ref{zylcodazzi}) follows directly from (\ref{zylweingarten}) and 
the Codazzi-Mainardi equation \cite[p.~115]{oneill83a}
\begin{eqnarray*}
\<R^\cZ(X,Y)U,\nu\> &=& \<(\nMt_XW)(Y),U\> - \<(\nMt_YW)(X),U\> .
\end{eqnarray*}

The Riccati equation for $W$
$$
(\nZ_\nu W)(X) = R^\cZ(X,\nu)\nu + W^2(X)
$$
gives
\begin{eqnarray*}
\<R^\cZ(X,\nu)\nu,Y\> &=&
\<(\nZ_\nu W)(X),Y\> - \<W^2(X),Y\> \\
&=&
\pa{}{t}\<W(X),Y\> - \<W(\nZ_\nu X),Y\> - \<W(X),\nZ_\nu Y\> \\
&&
+\half \gdt(W(X),Y) \\
&=&
-\half\pa{}{t}\gdt(X,Y) - \<W(\nZ_X \nu),Y\> - \<W(X),\nZ_Y \nu\> \\
&&
+\half \gdt(W(X),Y) \\
&=&
-\half\gddt(X,Y) + \<W(W(X)),Y\> + \<W(X),W(Y)\> \\
&&
+\half \gdt(W(X),Y) \\
&=&
-\half\gddt(X,Y) -\half \gdt(W(X),Y)
\end{eqnarray*}
which is (\ref{zylriccati}).
The Ricci curvature is now easily computed.
\begin{eqnarray*}
\ricZ(\nu,\nu) &=&
\summe \eps_i \<R^\cZ(e_i,\nu)\nu,e_i\> \\
&\stackrel{(\ref{zylriccati})}{=}&
-\half \summe \eps_i \left( \gddt(e_i,e_i) + \gdt(W(e_i),e_i) \right) \\
&\stackrel{(\ref{zylweingarten})}{=}&
-\half \tr_{g_t}(\gddt) + \tr(W^2)
\end{eqnarray*}
which is (\ref{zylricnn}).
Moreover,
\begin{eqnarray*}
\ricZ(X,\nu) &=&
\summe \eps_i \<R^\cZ(X,e_i)e_i,\nu\> \\
&\stackrel{(\ref{zylcodazzi})}{=}&
\half\summe\eps_i \left((\nMt_{e_i}\gdt)(X,e_i) 
- (\nMt_X\gdt)(e_i,e_i)\right)\\
&\stackrel{(\ref{zylweingarten})}{=}&
-\summe\eps_i \left(\<(\nMt_{e_i}W)(X),e_i\> 
- \<(\nMt_XW)(e_i),e_i\>\right)\\
&=&
-\<\divMt ,X\> + \tr(\nMt_X W)\\
&=&
-\<\divMt ,X\> + d_X\tr(W)
\end{eqnarray*}
thus showing (\ref{zylrictn}).
Furthermore,
\begin{eqnarray*}
\ricZ(X,Y)
&=&
\summe\eps_i \<R^\cZ(e_i,X)Y,e_i\> + \<R^\cZ(\nu,X)Y,\nu\> \\
&\stackrel{(\ref{zylgauss}),(\ref{zylriccati})}{=}&
\summe\eps_i \Big(\<R^{M_t}(e_i,X)Y,e_i\>
+ \frac{1}{4}\gdt(e_i,Y)\gdt(X,e_i)\\
&&
-\frac{1}{4}\gdt(e_i,e_i) \gdt(X,Y)\Big)
-\half\left( \gddt(X,Y) + \gdt(W(X),Y)  \right) \\
&=&
\ricMt(X,Y) + \summe\eps_i (\<W(e_i),Y\>\<W(X),e_i\> \\
&&
- \<W(e_i),e_i\> \<W(X),Y\>) -\half \gddt(X,Y) + \<W^2(X),Y\> \\
&=&
\ricMt(X,Y) + 2 \<W(X),W(Y)\> -\tr(W)\<W(X),Y\> \\
&&
-\half \gddt(X,Y)
\end{eqnarray*}
shows (\ref{zylrictt}).
Formula (\ref{zylscal}) for the scalar curvature follows from
(\ref{zylricnn}) and (\ref{zylrictt}).
\end{proof}

\Example{
A simple special case of a generalized cylinder is that of a 
{\em warped product}, i.~e.\ $g_t = f(t)^2g$ where $f:I\to\RR$ is
a smooth positive function.
Then $\gdt = 2\,f\,\dot{f}\,g = \frac{2\fd}{f}g_t$ and 
$\gddt = 2(\fd^2 + f\fdd)g = 2\frac{\fd^2+f\fdd}{f^2}g_t$ and the
formulas in Proposition~\ref{zylformeln} reduce to
\begin{eqnarray*}
W &=& -\frac{\fd}{f} \id, \\
\RZ(X,Y)U &=& 
R^{M_t}(X,Y)U + \frac{\fd^2}{f^2}\left(\<X,U\>Y - \<Y,U\>X\right), \\
\RZ(X,\nu)\nu &=& - \frac{\fdd}{f} X, \\
\ricZ(X,Y) &=& \ricMt(X,Y) - \frac{(n-1)\fd^2 + f\fdd}{f^2}\<X,Y\>,\\
\ricZ(X,\nu) &=& 0,\\
\ricZ(\nu,\nu) &=& - n \frac{\fdd}{f},\\
\scalZ &=& \scal^{M_t} -\,\, n\, \frac{(n-1)\fd^2 + 2f\fdd}{f^2},
\end{eqnarray*}
compare \cite[Ch.~7]{oneill83a}.
}

%%%%%%%%%%%%%%%%%%%%%%%%%%%%%%%%%%%%%%%%%%%%%%%%%%%%%%%%%%%%%%%%%%%%%%%%%
\section{Identifying spinors and the variation formula for the Dirac operator}
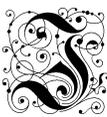
\label{secident}
%%%%%%%%%%%%%%%%%%%%%%%%%%%%%%%%%%%%%%%%%%%%%%%%%%%%%%%%%%%%%%%%%%%%%%%%%

\yinipar{I}t is an annoying problem that the definition of spinors, in contrast
to that of differential forms and tensors, depends on the semi-Riemannian
metric of the manifold.
Hence if one wants to compare the Dirac operators for two different
metrics one first has to identify the underlying spinor bundles.

The problem of constructing such identifications can be split into two
steps:
First construct identifications for any two metrics in a 1-parameter 
family of metrics.
The identification of spinors for two metrics will in general depend
on the 1-parameter family of metrics joining them.
Secondly, given two metrics construct a natural curve of metrics joining
them.

Both steps have been carried out very satisfactorily for the case of
Riemannian metrics in \cite{bourguignon-gauduchon92a}.
In the present section we will deal only with the first step.
The second step cannot always be carried out.
In Section~\ref{seclometrics} we will discuss this 
problem for the case of Lorentz metrics in great detail.

Now let $g_t$, $t\in I$, be a smooth 1-parameter family of semi-Riemannian
metrics of signature $(r,s)$ on a manifold $M$.
We form the generalized cylinder $\cZ := I \times M$ with metric 
$g = dt^2 + g_t$.
For $t\in I$ we abbreviate the semi-Riemannian manifold $(M,g_t)$ by $\Mt$.

%For any $t_0,t_1\in I$ parallel translation along the curves $t \mapsto
%(t,x) \in \cZ$, $x \in M$ fixed, provides an $\SOrs$-equivariant 
%isomorphism $\tau_{t_0}^{t_1} : \PSO(\Mto) \to \PSO(\Mtl)$.
%Given a spin structure $\Theta_{t_0} : \PSpin(\Mto) \rightarrow \PSO(\Mto)$ 
%on $\Mto$ we define a spin structure on $\Mtl$ by
%$$
%\Theta_{t_1} := \tau_{t_0}^{t_1} \circ \Theta_{t_0} :
%\PSpin(\Mtl) := \PSpin(\Mto) \rightarrow \PSO(\Mtl).
%$$
%This way we have constructed a smooth $\Spinrs$-principal fiber bundle
%$\Pt:=\bigcup_{t\in I}\PSpin(\Mt)$ over $\cZ$ and an equivariant map from 
%$\Pt$ to the bundle $P := \bigcup_{t\in I}\PSO(\Mt)$ of frames on $\cZ$ 
%tangent to the leaves $\Mt$.
%This map restricts to the spin structure $\Theta_{t} : \PSpin(\Mt) 
%\rightarrow \PSO(\Mt)$ for each leaf $\Mt$.
%Now $P$ injects naturally into $\PSO(\cZ)$ by $(e_1,\ldots,e_n) \mapsto
%(\nu,e_1,\ldots,e_n)$ where $\nu = \pa{}{t}$.
%Extending the structure group $\Spinrs$ of $\Pt$ to $\Spin(r+1,s)$ yields
%a spin structure $\PSpin(\cZ)$ on $\cZ$ which induces the spin structure
%$\PSpin(\Mt)$ on each leaf $\Mt$.

Spin structures on $M$ and on $\cZ$ are in 1-1-correspondence.
As explained in Section~\ref{secfoliation} spin structures on $\cZ$ can be 
restricted to spin structures on $\Mt = M$.
Conversely, given a spin structure on $M$ it can be pulled back
to $I \times M$ yielding a $\GLt^+(n,\RR)$-principal bundle on $\cZ$.
Enlarging the structure group via the embedding $\GLt^+(n,\RR) \hookrightarrow
\GLt^+(n+1,\RR)$ covering the standard embedding $\GL^+(n,\RR) \hookrightarrow
\GL^+(n+1,\RR)$, $a \mapsto \begin{pmatrix} 1 & 0 \\ 0 & a \end{pmatrix}$,
yields the spin structure on $\cZ$ which restricts to the given spin
structure on $M$.

Let us write ``$\cdot$'' for the Clifford multiplication on $\cZ$ and
``$\cdott$'' for the Clifford multiplication on $\Mt$.
Recall from Section~\ref{secfoliation} that $\Sigma\cZ|_{\Mt} = 
\Sigma\Mt$ as Hermitian vector bundles if $n=r+s$ is even and 
$\Sigma^+\cZ|_\Mt = \Sigma\Mt$ if $n$ is odd.
In both cases the Clifford multiplications are related by 
$X \cdott \phi = \nu \cdot X \cdot \phi$.
For given $x\in M$ and $t_0,t_1 \in I$ parallel translation on $\cZ$
along the curve $t \mapsto (t,x)$ is a linear isometry $\ttt : \Sigma_x\Mto
\to \Sigma_x\Mtl$.
Since ``$\cdot$'' and $\nu$ are parallel along the curve $t \mapsto (t,x)$ 
so is the family of Clifford multiplications ``$\cdott$'' and $\ttt$
preserves Clifford multiplication in the following sense:
$$
\ttt(X\bullet_{t_0}\phi) = (\ttt X)\bullet_{t_1}(\ttt\phi) .
$$
In general, the covariant derivative and hence parallel transport
depends on the semi-Riemannian metric and its first derivatives.
We note here that for fixed $x \in M$ the parallel transport $\ttt : 
T_x\Mto \to T_x\Mtl$ or $\ttt : \Sigma_x\Mto \to \Sigma_x\Mtl$ is determined 
by $g_t(x)$ and $\gdt(x)$, no $x$-derivatives of $g_t$ enter.
Namely, if $x^1, \ldots,x^n$ are local coordinates on $M$ and 
$X(t,x) = \sum_{j=1}^n\xi^j(x,t)\pa{}{x^j}$
is a parallel vector field along $t \mapsto (t,x)$, then this means
by (\ref{Gammazyl0}) and (\ref{zylweingarten})
\begin{eqnarray*}
0 &=& \frac{\nabla}{dt} X 
\,\,=\,\, \sum_{j=1}^n\left( \dot{\xi}^j 
+ \sum_{k=1}^n\Gamma_{k,0}^{\cZ,j} \xi^k \right)\pa{}{x^j} \\
&=&
\sum_{j=1}^n\left( \dot{\xi}^j 
+ \frac{1}{2}\sum_{k,\ell=1}^n g_t^{j\ell}{\dot{g}}_{t,k\ell}\xi^k \right)
\pa{}{x^j}.
\end{eqnarray*}
Thus $\ttt : T_x\Mto \to T_x\Mtl$ is given by solving the system of 
ordinary differential equations
$$
\dot{\xi}^j(t,x) = 
- \frac{1}{2}\sum_{k,\ell=1}^n g_t^{j\ell}(x){\dot{g}}_{t,k\ell}(x)\xi^k(t,x).
$$
For spinors the situation is similar.
By \cite[Prop.~2]{bourguignon-gauduchon92a} this shows that our identification
$\ttt$ of spinors for different metrics coincides with the one in 
\cite{bourguignon-gauduchon92a}.

Now we rewrite the commutator formula of Proposition~\ref{kommutator}.
For a section $\phi$ of $\Sigma\cZ$ (or $\Sigma^+\cZ$ if $n$ is odd)
we have 
\begin{equation}
[\nSZ_\nu,\DDMt]\,\phi = 
\DWt\phi - \frac{n}{2}\, \gradMt(H_t)\cdott\phi 
+ \half\,\divMt(W_t)\cdott\phi
\label{kommneu}
\end{equation}
where $\DDMt$ is the Dirac operator of $\Mt$, $\gradMt$ is the gradient
and $\divMt$ the divergence (of endomorphisms) on $\Mt$, $W_t$ is the
Weingarten map of $\Mt$ in $\cZ$ and $H_t = \frac1n \tr(W_t)$ the mean
curvature and finally $\DWt\phi = \sum_{i=1}^n\eps_i\, e_i\cdott
\nabla^{\Sigma\Mt}_{W_t(e_i)}\phi$ for any orthonormal basis $e_1,\ldots,e_n$.
From (\ref{zylweingarten}) we have $\divMt(W_t) = -\half \divMt(\gdt)$,
$H_t = -\frac{1}{2n} \tr_{g_t}(\gdt)$ and $\DWt = -\half \Dgdt$ where
$\Dgdt\phi = \sum_{i,j=1}^n \eps_i\eps_j \gdt(e_i,e_j)e_i\cdott
\nabla_{e_j}^{\Sigma\Mt}\phi$.
Thus (\ref{kommneu}) can be rewritten as 
\begin{equation}
[\nSZ_\nu,\DDMt]\,\phi = 
-\half\Dgdt\phi + \frac{1}{4}\, \gradMt(\tr_{g_t}(\gdt))\cdott\phi 
- \frac14\,\divMt(\gdt)\cdott\phi .
\label{kommneuneu}
\end{equation}
Now if $\phi$ is parallel along the curves $t \mapsto (t,x)$, i.~e.\ it
is of the form $\phi(t,x) = \tau_{t_0}^t \psi(x)$ for some spinor field
$\psi$ on $\Mto$, then the left hand side of (\ref{kommneuneu}) is at $t=t_0$
\begin{eqnarray*}
[\nSZ_\nu,\DDMt]\,\phi &=& \nSZ_\nu\DDMt\,\phi 
= \left.\frac{d}{dt}\right|_{t=t_0} \tau_t^{t_0}\DDMt\,\phi \\
&=&
\left.\frac{d}{dt}\right|_{t=t_0} \tau_t^{t_0}\DDMt \tau_{t_0}^t \psi .
\end{eqnarray*}
We have shown the variation formula for the Dirac operator:
\begin{thm}
\label{Diracvariation}
Let $g_t$ be a smooth 1-parameter family of semi-Riemannian metrics on a 
spin manifold $M$.
We write briefly $\Mt$ for the semi-Riemannian spin manifold $(M,g_t)$.
Let $\ttt$ be the identification of spinor spaces for $\Mto$ and $\Mtl$ 
constructed above,
let $\DDMt$ be the Dirac operator of $\Mt$, let ``$\cdott$'' be Clifford
multiplication on $\Mt$ and let $\Dgdt\phi = \sum_{i,j=1}^n \eps_i\eps_j 
\gdt(e_i,e_j)e_i\cdott\nabla_{e_j}^{\Sigma\Mt}\phi$.

Then for any smooth spinor field $\psi$ on $\Mto$ we have
$$
\left.\frac{d}{dt}\right|_{t=t_0} \tau_t^{t_0}\DDMt \tau_{t_0}^t \psi
=
-\half\mathfrak{D}^{\dot{g}_{t_0}}\psi + \frac{1}{4}\, \grad^{\Mto}
(\tr_{g_{t_0}}(\dot{g}_{t_0}))\bullet_{t_0}\psi 
- \frac14\,\div^{\Mto}(\dot{g}_{t_0})\bullet_{t_0}\psi .
$$
\end{thm}
This is exactly the formula given in \cite[Thm.~21]{bourguignon-gauduchon92a}
for Riemannian manifolds.

%%%%%%%%%%%%%%%%%%%%%%%%%%%%%%%%%%%%%%%%%%%%%%%%%%%%%%%%%%%%%%%%%%%%%%%%%
\section{Energy-momentum tensors}
\label{secenergymomentum}
%%%%%%%%%%%%%%%%%%%%%%%%%%%%%%%%%%%%%%%%%%%%%%%%%%%%%%%%%%%%%%%%%%%%%%%%%

\yinipar{T}heorem~\ref{Diracvariation} can be used to compute the
energy-momentum tensor for spinors.
In order to explain what this means we briefly sketch Lagrangian
field theory, see \cite[p.~153 ff]{deligne99a} for a more detailed 
introduction.
Let $M$ denote a differentiable manifold and let $\cG$ be a set of (smooth)
semi-Riemannian metrics on $M$, open in the $C^\infty$-topology.
Let $\pi : E \to \cG \times M$ be a fiber bundle with finite dimensional fibers.
For example, if $M$ carries a spin structure the fiber over $(g,x)\in 
\cG\times M$ could be the spinor space at $x$ with respect to the metric
$g$, $E_{(g,x)} = \Sigma_x^gM$.
For each fixed $g\in\cG$ the restriction $\pi\hml(\{g\}\times M) \to M$
is a fiber bundle over $M$ and we can form the space of smooth sections
$\cS_g$ of this bundle.
These Fr{\'e}chet manifolds $\cS_g$ give rise to a Fr{\'e}chet fiber bundle
$\cS := \bigcup_{g\in\cG}\cS_g \to \cG$.
Let $\cF\subset \cS$ be a Fr{\'e}chet submanifold such that the restriction
$\pi : \cF \to \cG$ is again a Fr{\'e}chet fiber bundle.

Now let $L:\cF \to \dichtM$ be a smooth map where $\dichtM$ denotes
the space of smooth densities on $M$, i.~e.\ smooth sections of $\Lambda^nT^*M
\otimes \mathfrak{o}_M$ where $\mathfrak{o}_M$ is the orientation line bundle.
We assume that $L$ is local in the sense that for $\phi\in\cF$ the density
$L(\phi)$ evaluated at $x\in M$ depends only on $\phi(x)$ and the 
$M$-derivatives of $\phi$ at $x$.
In other words, $L(\phi)(x)$ is a function of the jet $j^\infty_M\phi(x)$.
We call $L$ the {\em Lagrangian density}.
In physics it is customary to integrate over $M$ and call $\int_M L(\phi)$
the {\em Lagrangian} or the {\em action}.
We avoid this integration since in general the integral $\int_M L(\phi)$
need not exist.

We call a smooth 1-parameter family $\phi_t\in\cF_g$ with $\phi_0=\phi$ 
{\em compactly supported} if it is constant outside a compact subset 
$K\subset M$, i.~e.\ $\phi_t(x)=\phi(x)$ for all $x\in M \setminus K$
and all $t$.
Since $L$ is local $L(\phi_t)$ is constant outside $K$ as well so that
$\int_M (L(\phi_t) - L(\phi))$ exists and
$$
\left.\frac{d}{dt}\right|_{t=0}\int_M (L(\phi_t) - L(\phi)) 
=
\int_M\left.\frac{d}{dt}\right|_{t=0}L(\phi_t) .
$$
The section $\phi\in\cF_g$ is called {\em critical} for $L$ if for each
compactly supported deformation $\phi_t$
$$
\int_M\left.\frac{d}{dt}\right|_{t=0}L(\phi_t) = 0 .
$$

To explain the concept of energy-momentum tensors we need one more piece of
structure.
Let $H \subset T\cF$ be a connection.
This means that for any $\phi\in\cF$ we have 
$T_\phi\cF = T_\phi(\cF_{\pi(\phi)}) \oplus H_\phi$ and the restriction 
$d\pi|_{H_\phi} : H_\phi \to T_{\pi(\phi)}\cG$ is an isomorphism.
For fixed $\phi\in\cF$ and $g:=\pi(\phi)$ we have the linear map
$dL \circ (d\pi|_{H_\phi})\hml : T_g\cG \to \dichtM$.
Recall that $T_g\cG$ is nothing but the space of smooth $(2,0)$-tensors.
A smooth symmetric $(2,0)$-tensor $Q_\phi$ will be called the 
{\em energy-momentum tensor} for $\phi$ with respect to the Lagrangian $L$ if 
$$
dL \circ (d\pi|_{H_\phi})\hml(k) = \<Q_\phi,k\>_g\, dV_g
$$
for all $k\in T_g\cG$.
Here $\<\cdot,\cdot\>_g$ denotes the (pointwise) metric on symmetric 
$(2,0)$-tensors induced by $g$ and $dV_g$ is the Riemannian volume 
measure for $g$.
If it exists $Q_\phi$ is obviously unique.
By its definition the energy-momentum tensor describes the behavior of
the Lagrangian under variations of the metric.

\Example{
Let $M$ carry a spin structure, let $\cG$ be the set of all semi-Riemannian
metrics on $M$ and let $E$ be the universal spinor bundle, $E_{(g,x)} =
\Sigma_x^gM$.
Then $\cS$ is the universal bundle of spinor fields and we put $\cF := \cS$.
We fix $\lambda\in\RR$ and we define the Lagrangian $L$ by
$$
L(\phi) := \Re\<\phi,(D^g-\lambda)\phi\>_g \, dV_g
$$
where $D^g$ is the Dirac operator with respect to the metric $g=\pi(\phi)$.
If $\phi_t$ is a compactly supported deformation of $\phi$ we write
$\frac{d}{dt}|_{t=0}\phi_t=\phid$ and we compute
\begin{eqnarray*} 
\int_M\left.\frac{d}{dt}\right|_{t=0}L(\phi_t) 
&=&
\int_M \Re (\<\phid,(D^g-\lambda)\phi\>_g + 
\<\phi,(D^g-\lambda)\phid\>_g)\, dV_g \\
&=&
2\Re \int_M \<\phid,(D^g-\lambda)\phi\>_g\, dV_g .
\end{eqnarray*}
Thus $\phi$ is critical if and only if $(D^g-\lambda)\phi = 0$, i.~e.\
if $\phi$ is a Dirac-eigenspinor for the eigenvalue $\lambda$.

The connection $H$ is determined by the parallel translation $\ttt$
used in the previous section to identify spinors for different metrics.
More precisely, $H_\phi$ is the set of all $\left.\frac{d}{dt}\right|_{t=0}
\tau_0^t\phi$ for all smooth curves $g_t$ of metrics
with $g_0 = \pi(\phi)$.

Now let $g_t$ be such a 1-parameter family of metrics and write 
$k := \dot{g}_0$.
We compute
\begin{eqnarray*}
\lefteqn{dL \circ (d\pi|_{H_\phi})\hml(k) } \\
&=&
\left.\frac{d}{dt}\right|_{t=0} L(\tau_0^t\phi) \\
&=&
\left.\frac{d}{dt}\right|_{t=0} 
\Re\<\tau_0^t\phi,(D^{g_t}-\lambda)(\tau_0^t\phi)\>_{g_t} \, dV_{g_t}\\
&=&
\left.\frac{d}{dt}\right|_{t=0} 
\Re\<\phi,(\tau_t^0D^{g_t}\tau_0^t-\lambda)\phi\>_{g_0}\, 
\frac{dV_{g_t}}{dV_{g_0}} \, dV_{g_0} \\
&=&
\Re \left(\<\phi,\left.\frac{d}{dt}\right|_{t=0}
(\tau_t^0D^{g_t}\tau_0^t\phi)\>_{g_0}
+ \<\phi,(D^{g_0}-\lambda)\phi\>_{g_0}
\left.\frac{d}{dt}\right|_{t=0}\frac{dV_{g_t}}{dV_{g_0}}\right) 
\, dV_{g_0} .
\end{eqnarray*}
The first term is given by the variation formula for the Dirac operator.
Since Clifford multiplication with tangent vectors is skewadjoint all
terms of the form $\Re\<\phi,X\bullet_{g_0}\phi\>$ vanishes.
Thus Theorem~\ref{Diracvariation} yields
$$
\Re \<\phi,\left.\frac{d}{dt}\right|_{t=0}(\tau_t^0D^{g_t}\tau_0^t\phi)\>_{g_0}
=
-\half\Re\<\phi,\mathfrak{D}^{k}\phi\>_{g_0} .
$$
For the second term we use
$$
\left.\frac{d}{dt}\right|_{t=0}\frac{dV_{g_t}}{dV_{g_0}} = \half \tr_{g_0}(k) .
$$
Thus 
\begin{eqnarray*}
dL \circ (d\pi|_{H_\phi})\hml(k)
&=&
\half\Re\left(-\<\phi,\mathfrak{D}^{k}\phi\>_{g_0} 
+ \<\phi,(D^{g_0}-\lambda)\phi\>_{g_0}\tr_{g_0}(k) \right) \, dV_{g_0} \\
&=&
\<Q_\phi,k\>_{g_0}\, dV_{g_0}
\end{eqnarray*}
for the symmetric $(2,0)$-tensor
\begin{eqnarray*}
Q_\phi(X,Y) &=& 
-\frac14 \Re\left(\<\phi,X\bullet_{g_0}\nabla_Y^{\Sigma M}\phi\>
+ \<\phi,Y\bullet_{g_0}\nabla_X^{\Sigma M}\phi\>\right) \\
&&
+\half\Re\<\phi,(D^{g_0}-\lambda)\phi\>\, g_0(X,Y) .
\end{eqnarray*}
If $\phi$ is critical, i.~e.\ if $D^{g_0}\phi = \lambda \phi$, then the 
energy-momentum tensor simplifies to
\begin{equation}\label{qphi}
Q_\phi(X,Y) = 
-\frac14 \Re\left(\<\phi,X\bullet_{g_0}\nabla_Y^{\Sigma M}\phi\>
+ \<\phi,Y\bullet_{g_0}\nabla_X^{\Sigma M}\phi\>\right) . 
\end{equation}
}

\Example{
Again, let $M$ carry a spin structure, let $\cG$ be the set of all 
semi-Riemannian metrics on $M$ and let $E$ be the universal spinor bundle, 
$E_{(g,x)} = \Sigma_x^gM$.
Then again $\cS$ is the universal bundle of spinor fields and we this time
we put $\cF_g := \{ \phi\in\cS_g\ |\ \int_M\<\phi,\phi\>_gdV_g=\pm 1 \}$.
We define the Lagrangian $L$ by
$$
L(\phi) := \Re\<\phi,D^g\phi\>_g \, dV_g .
$$
Now $\phi$ is critical if and only if
$$
\int_M\left.\frac{d}{dt}\right|_{t=0}L(\phi_t) =
2\Re \int_M \<\phid,D^g\phi\>_g\, dV_g = 0
$$
for all $\phid$ perpendicular to $\phi$, i.~e.\ if and only if $D^g\phi$
is a multiple of $\phi$.
This way we obtain all nonnull eigenspinors for all eigenvalues simultaneously
as critical $\phi$'s.

This time the connection has to be chosen differently because $\ttt$
is a pointwise isometry but the volume element $dV_g$ also depends on
the semi-Riemannian metric.
Therefore $\ttt$ does not give an isometry for the $L^2$-product used
to define $\cF$.
This can be corrected by defining the connection $\bar{H}$ as the set of 
all $\left.\frac{d}{dt}\right|_{t=0} \sqrt{\frac{dV_{g_t}}{dV_{g_0}}}
\tau_0^t\phi$ for all smooth curves $g_t$ of metrics with $g_0 = \pi(\phi)$.

Then we have for such a 1-parameter family of metrics $g_t$ with 
$k := \dot{g}_0$
$$
dL \circ (d\pi|_{\bar{H}_\phi})\hml(k)
=
\Re \<\phi,\left.\frac{d}{dt}\right|_{t=0}
(\tau_t^0D^{g_t}\tau_0^t\phi)\>_{g_0}
\, dV_{g_0} 
$$
and therefore
$$
Q_\phi(X,Y) = 
-\frac14 \Re\left(\<\phi,X\bullet_{g_0}\nabla_Y^{\Sigma M}\phi\>
+ \<\phi,Y\bullet_{g_0}\nabla_X^{\Sigma M}\phi\>\right)
$$
for all $\phi$, critical or not.
}
These two examples show that for noncritical $\phi$ the energy-momentum 
tensor also depends on the choice of connection $H$.
In contrast, for critical $\phi$ the differential $dL$ descends to 
a map $dL : T_\phi\cF/T_\phi(\cF_{\pi(\phi)}) \to \dichtM$.
Thus the map $dL \circ d\pi\hml : T_{\pi(\phi)}\cG \to \dichtM$ is 
well defined without any reference to $H$.

%%%%%%%%%%%%%%%%%%%%%%%%%%%%%%%%%%%%%%%%%%%%%%%%%%%%%%%%%%%%%%%%%%%%%%%%%
\section{Embeddings of hypersurfaces}
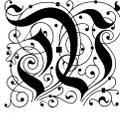
\label{sechyper}
%%%%%%%%%%%%%%%%%%%%%%%%%%%%%%%%%%%%%%%%%%%%%%%%%%%%%%%%%%%%%%%%%%%%%%%%%

\yinipar{W}e will now apply the cylinder construction described in 
Section~\ref{secgenzyl} to study the question whether a given manifold can 
be isometrically immersed as a hypersurface into a manifold of constant
curvature.
The classical example for such a result is the 
fundamental theorem for hypersurfaces which can be stated as follows:

\begin{thm} 
\label{fundthmhypersurface}
Let $(M^n,g)$ be a Riemannian manifold and let $A$ be a field 
of symmetric endomorphisms of $TM$ satisfying the equations of Gauss and 
Codazzi-Mainardi:
\begin{eqnarray}
(\nM_XA)Y &=& (\nM_YA)X,\label{aga1}\\
R^M(X,Y)Z &=&  \<A(Y),Z\>A(X) - \<A(X),Z\>A(Y)\label{aga2}
\end{eqnarray}
for all $X, Y, Z \in T_pM$, $p \in M$.

Then every point of $M$ has a neighborhood which can be isometrically 
embedded into Euclidean $(n+1)$-space 
$\RR^{n+1}$, with Weingarten map $A$. 
If $M$ is simply connected, then
there exists a global isometric immersion of $M$ into $\RR^{n+1}$ with the 
above property.
\end{thm}

A proof can be found in \cite[Ch.~VII.7]{kobayashi-nomizu69a}, but here 
we will give a more geometrical argument based on the cylinder construction. 
This will allow us to extend the result without effort to the semi-Riemannian
case and to embeddings into model spaces of constant sectional curvature not
necessarily zero.
We will construct an {\em explicit} metric of constant curvature on the 
cylinder $I\times M$, whose restriction to the leaf $\{0\}\times M$ is $g$. 

For a constant $\kappa\in\RR$ define the {\em generalized sine} and {\em 
cosine functions}
$$
  \asin_{\kappa}(t):=\left\{\begin{array}{c@{,\;\kappa\,}l}
      \frac{1}{\sqrt{\kappa}}\sin(\sqrt{\kappa}\cdot t)     & >0\\
      t                                                     & =0\\
      \frac{1}{\sqrt{|\kappa|}}\sinh(\sqrt{|\kappa|}\cdot t) & <0 
    \end{array}\right.
  \quad\mbox{ and }\quad
  \acos_{\kappa}(t):=\left\{\begin{array}{c@{,\;\kappa\,}l}
      \cos(\sqrt{\kappa}\cdot t)   & >0\\
      1                            & =0\\
      \cosh(\sqrt{|\kappa|}\cdot t) & <0 
    \end{array}\right.
$$
One easily checks $\asin_{\kappa}(0)=0$, $\acos_{\kappa}(0)=1$, 
$\kappa\asin_{\kappa}^2+\acos_{\kappa}^2=1$, $\asin_{\kappa}'=\acos_{\kappa}$,
and $\acos_{\kappa}'=-\kappa\asin_{\kappa}$.

\begin{thm}
\label{explicitembedding}
Let $(M^n,g)$ be a semi-Riemannian manifold and let $\kappa\in\RR$.
Let $A$ be a field of symmetric endomorphisms of $TM$ satisfying
\begin{eqnarray}
(\nM_XA)Y &=& (\nM_YA)X,\label{aga3}\\
R^M(X,Y)Z &=&  \<A(Y),Z\>A(X) - \<A(X),Z\>A(Y)\nonumber\\
&& + \,\kappa(\<Y,Z\>X-\<X,Z\>Y)
\label{aga4}
\end{eqnarray}
for all $X, Y, Z \in T_pM$, $p \in M$.
Define a family of metrics on $M$ by 
$$
g_t(X,Y) := g((\acos_\kappa(t)\,\id-\asin_\kappa(t)A)^2X,Y).
$$ 

Then the metric $dt^2+g_t$ on $\cZ=I\times M$ has constant sectional curvature 
$\kappa$ on its domain of definition (i.~e.\ for $|t|$ sufficiently small).
\end{thm}
\begin{proof}
Put $\RZk(X,Y)Z := \RZ(X,Y)Z - \,\kappa(\<Y,Z\>X-\<X,Z\>Y)$.
Having constant sectional curvature $\kappa$ is equivalent to $\RZk\equiv 0$.
The proof is based on the following lemma:

\begin{lemma} \label{a63} 
Let $\cZ=I\times M$ be a generalized cylinder and let $\kappa\in\RR$.
Assume that $g(\RZk(X,\nu)\nu,Y)=0$ for all vector fields $X$ and $Y$ 
on $\cZ$, where $\nu$ denotes the vector $\pa{}{t}$. 

{\rm (i)} 
If the Weingarten map $A$ of the hypersurface $\{0\}\times M$ of $\cZ$
satisfies (\ref{aga3}), then $g(\RZk(X,Y)Z,\nu)=0$ for all
vector fields $X$, $Y$ and $Z$ on $\cZ$.

{\rm (ii)}
If, moreover, $A$ also satisfies (\ref{aga4}), then $\RZk\equiv 0$, i.~e.\
$\cZ$ has constant sectional curvature $\kappa$.
\end{lemma}

Assume this lemma for a moment. 
We will check that the metric $dt^2+g_t$ satisfies the hypothesis of the 
lemma for $g_t(X,Y) = g((\acos_\kappa(t)\,\id-\asin_\kappa(t)A)^2X,Y)$. 
Let $W_t$ denote the Weingarten tensor of the hypersurface $\{t\}\times M$ 
of $\cZ$. 
This gives rise to a tensor field $W$ on $\cZ$, vanishing in the direction 
of $\nu$.
From the definition of $g_t$ we compute 
\begin{eqnarray*}
\gdt(X,Y) &=& -2 g((\acos_\kappa(t)\,\id-\asin_\kappa(t)A))
(\kappa\asin_\kappa(t)\,\id+\acos_\kappa(t)A)X,Y) \\
&=&
-2 g_t((\acos_\kappa(t)\,\id-\asin_\kappa(t)A))\hml
(\kappa\asin_\kappa(t)\,\id+\acos_\kappa(t)A)X,Y)
\end{eqnarray*}
hence by (\ref{zylweingarten})
$$
W = (\acos_\kappa(t)\,\id-\asin_\kappa(t)A))\hml
(\kappa\asin_\kappa(t)\,\id+\acos_\kappa(t)A).
$$
Moreover, 
\begin{eqnarray*}
\gddt(X,Y) &=&
-2 g\left([\kappa(\acos_\kappa(t)\,\id-\asin_\kappa(t)A)^2 
- (\kappa\asin_\kappa(t)\,\id+\acos_\kappa(t)A)^2]X,Y\right) .
\end{eqnarray*}
Equation (\ref{zylriccati}) yields
\begin{eqnarray*}
g_t(\RZ(X,\nu)\nu,Y) 
&=&
-\half \gddt(X,Y) -\half \gdt(W(X),Y) \\
&=&
g(\kappa(\acos_\kappa(t)\,\id-\asin_\kappa(t)A)^2X,Y) \\  
&=&
\kappa\, g_t(X,Y) ,
\end{eqnarray*}
thus $\RZ(X,\nu)\nu = \kappa\, X$ and hence $\RZk(X,\nu)\nu = 0$.
All conditions of the lemma are satisfied and the 
theorem follows. 
$\hfill\Box$

{\it Proof of the lemma.} 
The modified curvature tensor $\RZk$ has all the symmetries of a curvature
tensor including the Bianchi identities.

{\it i}) Consider the family of tensors on $M$ defined by
$K_t(X,Y,Z)_x:=\<\RZk(X,Y)Z,\nu\>_{(t,x)}$. Using the second Bianchi 
identity on $\cZ$, together with the fact that $\nu$ commutes 
with vectors on $M$ and the formula $W(X) = -\nZ_X \nu = -\nZ_\nu X + [\nu,X]
= -\nZ_\nu X$ we see
\begin{eqnarray}\dot{K_t}(X,Y,Z)&=&d_\nu\<\RZk(X,Y)Z,\nu\>\nonumber\\
&=&\<(\nZ_\nu\RZk)(X,Y)Z,\nu\>\nonumber\\
&&-\<\RZk(W(X),Y)Z+\RZk(X,W(Y))Z+\RZk(X,Y)W(Z),\nu\>\nonumber\\
&=&\<(\nZ_X\RZk)(\nu,Y)Z,\nu\>+\<(\nZ_Y\RZk)(X,\nu)Z,\nu\>\nonumber\\
&&+(W^*K_t)(X,Y,Z)
\label{akf}
\end{eqnarray}
where $W^*$ denotes the induced action of $W$ as a derivation on tensors.
From the assumption in the lemma we conclude
\begin{eqnarray*}
0 &=& d_X\<\RZk(\nu,Y)Z,\nu\> \\
&=&
\<(\nZ_X\RZk)(\nu,Y)Z,\nu\> + \<\RZk(\nZ_X\nu,Y)Z,\nu\> +
\<\RZk(\nu,\nZ_XY)Z,\nu\> \\
&&
+ \<\RZk(\nu,Y)\nZ_XZ,\nu\> + \<\RZk(\nu,Y)Z,\nZ_X\nu\> \\
&=&
\<(\nZ_X\RZk)(\nu,Y)Z,\nu\> - \<\RZk(W(X),Y)Z,\nu\> + 0 \\
&& +\, 0 - \<\RZk(\nu,Y)Z,W(X)\> 
\end{eqnarray*}
thus 
$$
\<(\nZ_X\RZk)(\nu,Y)Z,\nu\> = \<\RZk(W(X),Y)Z,\nu\> + \<\RZk(\nu,Y)Z,W(X)\>
$$
and similarly
$$
\<(\nZ_Y\RZk)(X,\nu)Z,\nu\> = \<\RZk(X,W(Y))Z,\nu\> + \<\RZk(X,\nu)Z,W(Y)\> .
$$
Plugging this into (\ref{akf}) yields
\begin{eqnarray*}
\dot{K_t}(X,Y,Z) &=& \<\RZk(W(X),Y)Z,\nu\>+\<\RZk(\nu,Y)Z,W(X)\> \\
&&+\<\RZk(X,W(Y))Z,\nu\>+\<\RZk(X,\nu)Z,W(Y)\>\\&&+(W^*K_t)(X,Y,Z) .
\end{eqnarray*}
Hence $\dot{K_t}=F(t)(K_t)$ for some linear endomorphism $F$ of the 
space of 3-tensors. 
This is a linear first order ODE for $K_t$.
The initial condition $K_0=0$ follows from 
(\ref{zylcodazzi}) because $W_0=A$ is a Codazzi tensor. This shows that
$K_t\equiv0$. 

{\it ii}) Similarly, using the identity $\<\RZk(X,Y)Z,\nu\> \equiv 0$ 
that we just obtained, we see that the family of tensors on $M$ defined by
$R_t(X,Y,Z,V)_x:=\<\RZk(X,Y)Z,V\>_{(t,x)}$ satisfies a linear ODE. Moreover,
(\ref{zylgauss}) implies $R_0\equiv 0$ because $W_0=A$ satisfies the 
Gauss equation. Thus $R_t\equiv 0$ for all $t$. This proves the lemma.
\end{proof}

Now recall that any semi-Riemannian manifold of constant sectional curvature 
$\kappa$ is locally isometric to $\MM_\kappa^{r,s}$.
Here $\MM_\kappa^{r,s}$ is the model space of constant sectional curvature
$\kappa$ and signature $(r,s)$.
If $\kappa=0$, then $\MM_0^{r,s}$ is semi-Euclidean space $\RR^n$ with the
metric $g_{r,s} = (dx^1)^2 + \cdots + (dx^r)^2 - (dx^{r+1})^2 - \cdots
- (dx^n)^2$.
If $\kappa>0$, then $\MM_\kappa^{r,s}$ is a pseudosphere, more precisely,
it is the semi-Riemannian hypersurface of $(\RR^{n+1},g_{r+1,s})$ defined
by $\<x,x\>_{r+1,s} = 1/\kappa$ and $x^1 > 0$ if $r=0$.
If $\kappa<0$, then  $\MM_\kappa^{r,s}$ is a pseudohyperbolic space, more 
precisely, it is the semi-Riemannian hypersurface of $(\RR^{n+1},g_{r,s+1})$ 
defined by $\<x,x\>_{r,s+1} = 1/\kappa$ and $x^{n+1} > 0$ if $r=0$.
In all cases $\MM_\kappa^{r,s}$ is connected and homogeneous.
Moreover, $\MM_\kappa^{r,s}$ is simpy connected except for 
$\MM_\kappa^{1,n-1}$ if $\kappa>0$ and $\MM_\kappa^{n-1,1}$ if $\kappa<0$, 
compare \cite[p.~108 ff]{oneill83a}.

The local isometry is essentially given by the Riemannian exponential map,
see \cite[Cor.~2.3.8]{wolf84a}, and it is uniquely determined by its
differential at a point.
Applying this to the cylinder constructed in 
Theorem~\ref{explicitembedding} yields the local statement in
the fundamental theorem for hypersurfaces for semi-Riemannian manifolds.

\begin{cor}
Let $(M^n,g)$ be a semi-Riemannian manifold of signature $(r,s)$ and
let $\kappa\in\RR$.
Let $A$ be a field of symmetric endomorphisms of $TM$ satisfying the 
equations of Gauss and Codazzi-Mainardi:
\begin{eqnarray*}
(\nM_XA)Y &=& (\nM_YA)X, \\
R^M(X,Y)Z &=& \<A(Y),Z\>A(X) - \<A(X),Z\>A(Y) \\
&& + \,\kappa(\<Y,Z\>X-\<X,Z\>Y)
\end{eqnarray*}
for all $X, Y, Z \in T_pM$, $p \in M$.

Then for every point $p\in M$, for every $q\in \MM_\kappa^{r+1,s}$, and for
every linear isometric embedding $F:T_pM \to T_q\MM_\kappa^{r+1,s}$ there
exists a neighborhood $U$ of $p$ in $M$ and an isometric embedding $f:U\to
\MM_\kappa^{r+1,s}$ as a semi-Riemannian hypersurface with Weingarten map $A$,
such that $f(p)=q$ and $df(p)=F$.

Moreover, any two such local embeddings $f_1$ and $f_2$ must agree
in a neighborhood of $p$ if $f_1(p)=f_2(p)=:q$ and $df_1(p)=df_2(p) :
T_pM \to T_q\MM_\kappa^{r+1,s}$.
\end{cor}

Now that this local result is established exactly the same proof as in
\cite[Ch.~VII, Thm.~7.2]{kobayashi-nomizu69a} can be used to show the
corresponding global immersion statement in the simply connected case.

\begin{cor}
Let $(M^n,g)$ be a simply connected semi-Riemannian manifold of signature 
$(r,s)$, let $\kappa\in\RR$ and let $A$ be a field of symmetric endomorphisms 
of $TM$ satisfying the two equations (\ref{aga3}) and (\ref{aga4}) above. 

Then $M$ can be isometrically immersed as a semi-Riemannian hypersurface 
into the model space $\MM_\kappa^{r+1,s}$ with Weingarten map $A$. 
Any two such immersions differ by an isometry of $\MM_\kappa^{r+1,s}$.
\end{cor}

%%%%%%%%%%%%%%%%%%%%%%%%%%%%%%%%%%%%%%%%%%%%%%%%%%%%%%%%%%%%%%%%%%%%%%%%%%
\section{Generalized Killing spinors}
\label{seckilling}
%%%%%%%%%%%%%%%%%%%%%%%%%%%%%%%%%%%%%%%%%%%%%%%%%%%%%%%%%%%%%%%%%%%%%%%%%%

\yinipar{W}e now turn our attention to restrictions of spinors to hypersurfaces. 
Let $M^n\subset \cZ^{n+1}$ be a hypersurface of a spin manifold $\cZ$ admitting 
a parallel 
spinor $\Psi$. If $n+1$ is even, we will assume that $\Psi$ lies in 
$\Sigma^+\cZ$. From the discussion in Section~\ref{secfoliation} 
we see that the restriction 
$\psi$ of $\Psi$ to $M$ is actually a spinor on $M$ and (\ref{spingauss}) 
reads
\begin{equation}\label{atk}
0=\n^\SZ_X\Psi=\n^\SM_X\psi-\frac{1}{2}A(X)\bullet\psi
\end{equation}
for all $X\in TM$ where $A$ is the Weingarten tensor of the submanifold $M$
and ``$\bullet$'' denotes Clifford multiplication on $M$. 
If $\psi$ is an eigenspinor of the Dirac operator, then $A$ is closely related
to the energy-momentum tensor of $\psi$.
More precisely, using (\ref{qphi}) one computes
$$
Q_\psi(X,Y) = \frac14 \<X,A(Y)\> \<\psi,\psi\>
$$
where $\<\psi,\psi\>$ is constant since $\psi$ is parallel on $\cZ$.
Spinors satisfying (\ref{atk}) will be called {\it generalized Killing spinors}. 
They are closely related to the so--called $T$--Killing spinors studied 
by Friedrich and Kim in \cite{friedrich-kim99a}.

Conversely, given a generalized Killing spinor $\psi$ on a manifold $M^n$ with 
$\n^\SM_X\psi-\frac{1}{2}A(X)\bullet\psi$, 
it is natural to ask whether the tensor $A$ can 
be realized as the Weingarten tensor of some isometric embedding of 
$M$ in a manifold $\cZ^{n+1}$ carrying parallel spinors. 
Morel studied this problem in the case where the tensor $A$ is parallel, see
\cite{morel03a}.

The next result provides an affirmative answer to the above question, 
for the case where the energy-momentum tensor of $\psi$ is a Codazzi tensor. 

\begin{thm} 
Let $(M^n,g)$ be a semi-Riemannian spin manifold and let $A$ be a field 
of symmetric endomorphisms of $TM$ satisfying equation (\ref{aga1}) on $M$. 
Let $\psi$ be a spinor on $(M^n,g)$ satisfying for all $X\in TM$
\begin{equation}
\n^\SM_X\psi=\frac{1}{2}A(X)\bullet\psi .
\label{atk2}
\end{equation}

Then the generalized cylinder $\cZ=I\times M$
with the metric $dt^2+g_t$, where $g_t(X,Y)=g((\id-tA)^2X,Y)$,
and with the spin structure inducing the given one on $\{0\}\times M$
by restriction has a parallel spinor, whose restriction to the leaf  
$\{0\}\times M$ is just $\psi$.
\end{thm}

\begin{proof} 
The spinor $\psi$ defines a spinor $\Psi$ on $\cZ$ by parallel 
transport along the geodesics $\RR\times \{x\}$. More precisely, we define 
$\Psi_{(0,x)}:=\psi_x$ via the identification 
$\Sigma_xM\cong\Sigma_{(0,x)}\cZ$ (resp. $\Sigma_{(0,x)}^+\cZ$ for $n$ odd) 
and $\Psi_{(t,x)}=\tau_0^t\Psi_{(0,x)}$. By construction we have 

\begin{equation}\label{ano}\n^\SZ_\nu\Psi\equiv 0\, \ \hbox{and}\ \ 
\n^\SZ_X\Psi|_{\{0\}\times M}= 0
\end{equation}
for all $X\in TM$.

The explicit form of the metrics $g_t$ yields $\<\RZ(X,\nu)\nu,Y\>=0$ 
on $\cZ$ for all $X$ and $Y$ tangent to $M$ as in the proof of 
Theorem~\ref{explicitembedding}.
Since the Codazzi equation (\ref{aga1}) holds Lemma~\ref{a63} (i) yields
$\<\RZ(\nu,X)Y,Z\>=0$ on all of $\cZ$.
Hence $R^\cZ(\nu,X)=0$ for all $X\in TM$. 

Let $X$ be a fixed 
arbitrary vector field on $M$, identified as usual with the vector field 
$(0,X)$ on $\cZ$. Using (\ref{ano}) we get 
$0=\half R^\cZ(\nu,X)\cdot\Psi=\n^\SZ_\nu\n^\SZ_X\Psi$, thus showing that 
the spinor field $\n^\SZ_X\Psi$ is parallel along the geodesics $\RR\times \{x\}$. 
Now (\ref{ano}) shows that this spinor vanishes for $t=0$, hence it is 
zero everywhere on $\cZ$. 
Since $X$ was arbitrary, this shows that $\Psi$ is parallel on $\cZ$.
\end{proof}  

This theorem generalizes the result from \cite{baer93a} where the case $A=\lambda
\cdot\id$ is treated, $\lambda\in\RR$, and it is shown 
that the cone over a manifold with Killing spinors admits parallel spinors, 
as well as a more recent result by Morel \cite{morel03a} for the case when 
$A$ is parallel.
Nevertheless, the question whether a manifold with a spinor satisfying 
(\ref{atk2}) can be isometrically embedded in a manifold with parallel spinors 
such that $A$ becomes the Weingarten tensor of the embedding without assuming
that $A$ is a Codazzi tensor is left open in the present article.

%%%%%%%%%%%%%%%%%%%%%%%%%%%%%%%%%%%%%%%%%%%%%%%%%%%%%%%%%%%%%%%%%%%%%%%%%
\section{The space of Lorentzian metrics}
\label{seclometrics}
%%%%%%%%%%%%%%%%%%%%%%%%%%%%%%%%%%%%%%%%%%%%%%%%%%%%%%%%%%%%%%%%%%%%%%%%%

\yinipar{I}n the final section we address the problem of connecting any two
semi-Riemannian metrics of signature 
$(r, s)$ on some manifold $M$ of dimension $n = r + s$,  by a curve $g
_t$ of semi-Riemannian metrics of the same signature in a unique and
universal manner.
The latter requirement reduces this problem to the
purely algebraic issue of finding a universal way of relating any two
inner products of signature 
$(r, s)$ on some real vector space $E \cong \mathbb{R}^{n}$ in
the manifold $\mathcal{M}_{r, s}$ of all inner products of
signature $(r, s)$ on $E$. 

In the positive or negative definite case an obvious candidate is the 
linear interpolation $g_t = t g_1 + (1 - t) g_0$
which, however, cannot be used for other signatures. An alternative
solution, which has been considered
in the definite case, see e.g. \cite{bourguignon-gauduchon92a}, but holds in 
a formally
identical way for all signatures, relies on the geometry of
$\mathcal{M}_{r, s}$, as a (semi-Riemannian) symmetric space that
we now recall briefly.  

For any signature $(r, s)$ the identity component of the general linear group 
$\GL^+  (E) \cong \GL^+ (n, \mathbb{R})$
acts transitively on $\mathcal{M}_{r, s}$ by
\begin{equation} 
(\gamma \cdot g) (u, v) = g (\gamma^{-1} u, \gamma
^{-1} v) \nonumber
\end{equation}
for $\gamma\in\GL^+ (E)$, $g\in\mathcal{M}_{r, s}$, and $u, v\in E$. 
For any chosen $g_0$ in $\mathcal{M}_{r, s}$, the
isotropy group of $g_0$ in $\GL^+ (E)$ is the special orthogonal
group $\SO(g_0)$ relative to $g_0$.
Recall that, except in  the definite case where $\SO(g_0)$ is connected, 
$\SO(g_0)$ has {\it two}  connected components. 
We thus get the identification $\mathcal{M}_{r, s} = \GL^+ (E)/\SO(g_0)$ 
or, equivalently, $\mathcal{M}_{r, s} = \mathbb{R}^+ \times \SL(E)/\SO(g_0)$,
where $\mathbb{R}^+$ acts by homotheties, and $\SL(E) \cong 
\SL(n, \mathbb{R})$ denotes the special linear group of elements of 
determinant $1$ in $\GL^+ (E)$.
Hence $\mathcal{M}_{r, s}^0 := \SL(E)/\SO(g_0)$ can be
regarded as the space of inner products on $E$ of signature $(r, s)$
and with a fixed volume element. 
Concerning the problem addressed in this section, it is clearly
sufficient to restrict our attention to $\mathcal{M}_{r, s}^0$. 

The homogeneous geometry of $\mathcal{M}_{r, s}^0 = \SL(E)/\SO(g_0)$ 
can be described as follows. 
For simplicity, write $G:=\SL(E)$, $H:=\SO(g_0)$, let $\mathfrak{g}$ be the 
Lie algebra of $G$, identified with the Lie algebra of trace-free 
endomorphisms of $E$, and let $\mathfrak{h}$ be the Lie algebra of $H$, 
identified with the Lie algebra of $g_0$-skewsymmetric endomorphisms.
Denote by $\mathfrak{m}$
the orthogonal complement of $\mathfrak{h}$ in $\mathfrak{g}$ with
respect to the Killing form of $\mathfrak{g}$, so that $\mathfrak{g} =
\mathfrak{h} \oplus \mathfrak{m}$. 
Recall that the Killing form of
$\mathfrak{g}$ equals the bilinear form $a, b \mapsto {\rm tr} (a
b)$, up to a positive universal constant, so that $\mathfrak{m}$ is
the space of $g_0$-symmetric elements of $\mathfrak{g}$. Since the
Killing form is $G$-invariant, $\mathfrak{m}$ is stable under the
adjoint action of $H$, making  $\mathcal{M}_{r, s}^0$ a
reductive homogeneous space. 
Moreover, we clearly have the Lie bracket relations $[\mathfrak{h},
\mathfrak{h}]\subset\mathfrak{h}$, $[\mathfrak{h},\mathfrak{m}]\subset
\mathfrak{m}$, and $[\mathfrak{m}, \mathfrak{m}]\subset\mathfrak{h}$ showing 
that $\mathcal{M}_{r, s}^0$ is actually a symmetric homogeneous space. 

In the positive definite case, $\mathcal{M}_{n, 0}^0$ is 
a Riemannian symmetric space of noncompact type, hence a
Hadamard space. It follows that any two points of $\mathcal{M}_{n, 0}
^0$  can be joined by a unique geodesic. If $g$ and $g_0$ are any two
points of $\mathcal{M}_{n, 0}$, then $g = g_0 (A \cdot, \cdot)$, for
a uniquely defined automorphism $A$ of $E$, where $A$ is symmetric and
positive definite for both $g_0$ and $g$. 
Then $A = {\rm exp} (a)$
for a uniquely defined symmetric endomorphism $a$ of $E$ and the unique
geodesic connecting $g_0$ to $g$ is the curve $g_t := g_0 ({\rm
  exp} (t a) \cdot, \cdot) = g_0 (A^t \cdot, \cdot)$, for $t\in [0, 1]$
where $\exp : \mathfrak{g} \to G$ denotes the exponential mapping. 

In the general case, the restriction of the Killing form to
$\mathfrak{m}$ is an $H$-invariant inner product of signature 
$\left(\frac{r (r + 1)}{2} + \frac{s (s + 1)}{2} - 1, r s\right)$, 
making $\mathcal{M}_{r, s}^0$  a {\it semi-Riemannian} symmetric space
of this signature.

The fact that $\mathcal{M}_{r, s}^0$ is symmetric, as a
semi-Riemannian homogeneous space, implies that the Levi-Civita
connection of the semi-Riemannian metric coincides with the
canonical homogeneous connection. 
In particular, all (semi-Riemannian) geodesics emanating from $g_0$ are of the form
\begin{equation} 
{\rm exp} (t X) \cdot g_0 \nonumber
\end{equation}
for $X\in\mathfrak{m} = T_{g_0} \mathcal{M}_{r, s}$.

As a homogeneous semi-Riemannian manifold $\mathcal{M}_{r, s}^0$ is
certainly geodesically complete in the sense that geodesics are defined on 
all of $\RR$, but for $(r, s) \neq (n, 0), (0, n)$,
it is not longer true that any two points can be joined by a geodesic
and, if so, there is no guarantee that the geodesic be unique. 
This will be illustrated firstly  in the case that 
$(r, s) = (1, 1)$, then in the general Lorentzian case when $(r, s) =
(n - 1, 1)$.  

%%%%%%%%%%%%%%%%%%%%%%%%%%%%%%%%%%%%%%%%%%%%%%%%%%%%%%%%%%%%%%%%%%%%%%%%%%%%%

\subsection{The space of Lorentzian inner products in dimension $2$} 
\label{ssl2}

%%%%%%%%%%%%%%%%%%%%%%%%%%%%%%%%%%%%%%%%%%%%%%%%%%%%%%%%%%%%%%%%%%%%%%%%%%%%%

Let $E$ denote an oriented real vector space of dimension $2$. 
We fix a positive generator $\omega$ of the real line $\Lambda^2 E^*$,
which can be viewed as a symplectic form on $E$.
Now $G \cong \SL(2,\mathbb{R})$, $\mathfrak{g} \cong \mathfrak{sl} 
(2, \mathbb{R})$ is the Lie algebra of trace-free endomorphisms of $E$, and 
$\mathcal{M}_{1, 1}^0$ is the space of all Lorentzian inner products on $E$, 
whose volume form with respect to the given orientation is $\omega$. 
For any chosen point $g_0 \in\mathcal{M}_{1, 1}^0$ 
we then have $\mathcal{M}_{1, 1}^0 = \SL(2,\RR)/\SO(1, 1)$. 
Note that $\SO(1,1)$ has two connected components. 
The connected component of the identity $\SO_0 (1,1)$ is isomorphic the the 
additive group $\mathbb{R}$ of real numbers via the isomorphism $t \mapsto
\begin{pmatrix} \cosh{t} & \sinh{t} \\ \sinh{t} & \cosh{t}\end{pmatrix}$.
The other connected component equals $- \SO_0 (1,1)$. 
Differentiation with respect to $t$ shows that the corresponding isotropy
Lie algebra $\mathfrak{h}$ is the Lie algebra of $2 \times
2$-matrices of the form $\begin{pmatrix} 0 & b \\ b & 0
\end{pmatrix}$, for $b\in\RR$. 

An endomorphism $\alpha$ of $E$ is tracefree if and only if it is
``antisymmetric'' with respect to $\omega$, i.~e.\ if and only if it
satisfies: $\omega (\alpha \cdot, \cdot) + \omega (\cdot, \alpha
\cdot) = 0$. 

For any $g \in \mathcal{M}_{1, 1}^0$ there is one and only one automorphism
${I_g}$ of $E$ such that
\begin{equation} 
\label{I} 
g = \omega (\cdot, {I_g} \cdot). 
\end{equation}

Since  $g$ is symmetric ${I_g}$ is trace-free.
Its determinant equals $-1$ because $g$ is Lorentzian, with volume form 
equal to $\omega$.
In particular, $I_g^2 = 1$. 
The light cone of $g$ is the union of the two eigenspaces of ${I_g}$, for
the eigenvalues $\pm 1$.
The latter are generated by $v \pm {I_g}v$ respectively, for any nonzero 
$v\in E$.

Conversely, for any automorphism $I$ of $E$ of trace equal to $0$ and
of determinant equal to $-1$, the bilinear form $g$ defined by
$g = \omega (\cdot, {I} \cdot)$ is a Lorentzian inner product, with 
volume form equal to $\omega$ and $I=I_g$.

The automorphism ${I_g}$ belongs to the Lie algebra $\mathfrak{g}$, 
on which $G$ acts by the adjoint representation, and the map 
$g \mapsto {I_g}$ is $G$-equivariant. 
Indeed, by definition of $G$, we have that $\omega (\gamma \cdot, \gamma \cdot)
= \omega (\cdot, \cdot)$ for each $\gamma \in G$, so that
\begin{equation} 
\gamma \cdot g 
= g (\gamma^{-1} \cdot, \gamma^{-1} \cdot) 
= \omega (\gamma^{-1} \cdot, {I_g} \, \gamma^{-1} \cdot) 
= \omega (\cdot, \gamma \, {I_g} \gamma^{-1} \cdot).
\nonumber\end{equation}

The map $g \mapsto {I_g}$ is then a $G$-equivariant identification of 
$\mathcal{M}_{1, 1}^0$ with the adjoint orbit of all elements 
of $\mathfrak{g}$ of determinant equal to $-1$.

As a function defined on $\mathfrak{g} \cong \RR^3$, the opposite of the 
determinant is a nondegenerate  quadratic form of signature (2, 1), equal 
to the (suitably normalized) Killing form. 
We denote the symmetric bilinear form corresponding to $-\det$ by 
$\<\cdot,\cdot\>$, i.~e.\ $\<u,u\> = - \det(u) = \half \tr(u^2)$.
The adjoint orbit is then the pseudosphere $\MM^{1,1}_1$ of elements $u$ 
such that $\<u,u\> = 1$ in the 3-dimensional Minkowski space $(\mathfrak{g},
\<\cdot,\cdot\> )$. 
The restriction of $\<\cdot,\cdot\>$ to $\MM^{1,1}_1$  makes the latter a 
$G$-homogeneous Lorentzian manifold, known as the 2-dimensional {\it de Sitter
  universe}.

The map $\mathcal{M}_{1, 1}^0 \to \MM^{1,1}_1$, $g \mapsto {I_g}$,
is a $G$-equivariant isometry.

Reflection with respect to $\<\cdot,\cdot\>$ about a vector subspace 
is an isometry of $(\mathfrak{g},\<\cdot,\cdot\> )$ and it preserves
$\MM^{1,1}_1$.
Since the fixed point set of an isometry is a totally geodesic submanifold
the geodesics of $\MM^{1,1}_1$ are precisely the intersections
of $\MM^{1,1}_1$ with 2-dimensional vector subspaces $E\subset\mathfrak{g}$.
There are three types of geodesics:
timelike geodesics (hyperbolas) corresponding to Minkowski planes,
spacelike geodesics (ellipses) corresponding to spacelike planes, 
and null geodesics (straight lines) corresponding to degenerate planes 
(tangent to the light cone).

\begin{center}
  \input{fig1c}
  \input{fig1a}
  \input{fig1b}
\end{center}
\begin{center}
{\em Fig.~1}  
\end{center}

Now let $I$, $I'$ be two different points in $\MM^{1,1}_1$.
If $I'=-I$, then each plane $E$ containing $I$ also contains $I'$.
In the timelike or in the null case $I'$ lies on the other connected
component of $E \cap \MM^{1,1}_1$.
Thus all spacelike geodesics emanating from $I$ hit $I'=-I$, but the
timelike and null geodesics emanating from $I$ miss $I'=-I$.

If $I'\neq -I$, then $I$ and $I'$ are linearly independent, so the plane $E$ 
containing $I$ and $I'$ is uniquely determined.
Thus $I'$ is hit by the geodesic emanating from $I$ if and only if 
it does not lie on the ``wrong'' connected component of $E \cap \MM^{1,1}_1$
(in the timelike or null case).
In other words, the points on $\MM^{1,1}_1$ which cannot be reached by a
geodesic emanating from $I$ are precisely the ones lying on timelike or
null geodesics emanating from $-I$.

\begin{center}
  \input{fig2}
\end{center}
\begin{center}
{\em Fig.~2}  
\end{center}

The two null geodesics emanating from $-I$ are cut out of $\MM^{1,1}_1$ by
the affine plane $\{\<I,I'\>=-1\}$.
Thus the points $I'\in\MM^{1,1}_1$ with $\<I,I'\> < -1$ cannot be attained
by a geodesic from $I$.

Similarly, by looking at the affine plane $\{\<I,I'\>=+1\}$ we see
that the points $I'$ with $\<I,I'\> > 1$ are the ones that lie on 
timelike geodesics emanating from $I$, the ones with $\<I,I'\> = 1$
are the ones that lie on null geodesics emanating from $I$, and 
the ones with $-1 < \<I,I'\> < 1$ lie on spacelike geodesics emanating 
from $I$.

We now retranslate this information back to $\mathcal{M}_{1, 1}^0$. 
If $g, g' \in\mathcal{M}_{1, 1}^0$, then
\begin{equation} g' = g (A \cdot, \cdot), \nonumber\end{equation}
with
\begin{equation} 
A = I_g^{-1} I_{g'} = I_g I_{g'}. \nonumber
\end{equation}
We then have
\begin{equation} 
\<I_{g}, I_{g'}\> = \frac{1}{2} {\rm tr}\, A. \nonumber
\end{equation}
Note that $A$ is $g$- and $g'$-symmetric and of determinant equal to
$+1$. 

By choosing $g$ as a base-point, we conclude that $\mathcal{M}_{1, 1}^0$ can 
also be identified with the space of all $g$-symmetric automorphisms of 
determinant $1$ of $E$. 
We summarize:

\begin{prop}
\label{prop2dim}
The space $\mathcal{M}_{1, 1}^0$ of Lorentzian inner products on a 2-dimensional
real vector space that have a fixed volume element carries a natural Lorentzian
metric making it $\SL(2,\RR)$-equivariantly isometric to the 2-dimensional
de Sitter universe.
For $g, g' \in\mathcal{M}_{1, 1}^0$ there is a unique endomorphism $A$ such that
$g' = g (A \cdot, \cdot)$.
Moreover, the following holds:
\begin{itemize}
\item
If $\tr(A) > 2$, then there is a unique geodesic in $\mathcal{M}_{1, 1}^0$ joining 
$g$ and $g'$.
This geodesic is timelike.
\item
If $\tr(A) = 2$, then there is a unique geodesic in $\mathcal{M}_{1, 1}^0$ joining 
$g$ and $g'$.
This geodesic is null.
\item
If $-2 <\tr(A) < 2$, then there is a unique geodesic in $\mathcal{M}_{1, 1}^0$ 
joining $g$ and $g'$.
This geodesic is spacelike.
\item
If $\tr(A) < -2$, then there is no geodesic in $\mathcal{M}_{1, 1}^0$ joining 
$g$ and $g'$.
\item
If $\tr(A) = -2$ and $g\not= -g'$, then there is no geodesic in 
$\mathcal{M}_{1, 1}^0$ joining $g$ and $g'$.
\item
If $\tr(A) = -2$ and $g = -g'$, then all spacelike geodesics in 
$\mathcal{M}_{1, 1}^0$ emanating from $g$ pass through $g'$ while the timelike 
and null geodesics in $\mathcal{M}_{1, 1}^0$ emanating from $g$ miss $g'$.
\end{itemize}
\end{prop}

This proposition shows that given two Lorentzian metrics on a 2-dimensional
manifold we can construct a canonical 1-parameter family of Lorentzian metrics
joining them only if the endomorphism field $A$ relating the two metrics
satisfies $\tr(A) > -2$.
A restriction like this does not come as a surprise because there are
pairs of Lorentzian metrics e.~g.\ on the 2-torus which cannot even be
joined by any continuous curve of Lorentzian metrics.

%%%%%%%%%%%%%%%%%%%%%%%%%%%%%%%%%%%%%%%%%%%%%%%%%%%%%%%%%%%%%%%%%%%%%%%%%%%%%%

\subsection{The space of Lorentzian inner products in higher dimensions}
\label{lorentzhighdim}

%%%%%%%%%%%%%%%%%%%%%%%%%%%%%%%%%%%%%%%%%%%%%%%%%%%%%%%%%%%%%%%%%%%%%%%%%%%%%%

We now consider the manifold $\mathcal{M}_{n - 1, 1} = \mathbb{R}^+
\times \mathcal{M}_{n - 1, 1}^0$ of all Lorentzian inner
products of signature $(n-1,1)$ on some $n$-dimensional real
vector space $E$. 

As observed before the manifold $\mathcal{M}_{n-1,1}^0$
is a symmetric semi-Riemannian space of signature $\left(\frac{n (n - 1)}{2},
n - 1\right)$ and the geodesics emanating from any chosen base-point $g_0$ are
of the form ${\rm exp} (t X) \cdot g_0$, where $X$ belongs to the
space $\mathfrak{m}$ of trace-free $g_0$-symmetric endomorphisms of
$E$, $\mathfrak{m}$ being naturally identified with the tangent space $T_{g_0} 
\mathcal{M}_{n-1,1}^0$.

The goal of this section is to determine the set of elements $g\in
\mathcal{M}_{n-1, 1}$ which can be joined from $g_0$ by a geodesic in 
$\mathcal{M}_{n-1,1}$, and whether or not this geodesic is unique. 
This has just been done in detail in the case that $n = 2$
and, as we shall  see, the general case can essentially be reduced to the
$2$-dimensional case.
More precisely, we have

\begin{prop} 
\label{proplorentz} 
Let $g_0$ and $g $ be two distinct points in $\mathcal{M}_{n - 1, 1}$. 
Then there is the following alternative: Either

{\rm (i)} 
$E$ splits as
\begin{equation} 
E = E_{1, 1} \oplus E_{n - 2, 0}, \nonumber
\end{equation}
where the sum is orthogonal,  $E_{1, 1}$ is of signature $(1, 1)$, $E_{n - 2, 0}$ 
is of signature $(n - 2, 0)$ for $g_0$ and $g$.
Both $g_0$ and $g$ belong to the corresponding totally geodesic submanifold 
$\mathcal{M}_{1, 1} \times \mathcal{M}_{n - 2, 0} \subset \mathcal{M}_{n-1, 1}$. 
Thus the issue of the existence and uniqueness of geodesics connecting $g_0$ to $g$ 
is reduced to the same issue for the $2$-dimensional Lorentzian
metrics ${g_0}_{|E_{1, 1}}$ and $g_{|E_{1, 1}}$ in $\mathcal{M}_{1, 1}$
as described in Proposition~\ref{prop2dim}, or 

{\rm (ii)} 
$E$ splits as
\begin{equation} 
E = E_{2, 1} \oplus E_{n - 3, 0}, \nonumber
\end{equation}
where the sum is orthogonal,  $E_{2, 1}$ is of signature $(2, 1)$, $E_{n - 3, 0}$ 
is of signature $(n - 3, 0)$ for $g_0$ and $g$.
Both $g_0$ and $g$ belong to the corresponding totally geodesic submanifold 
$\mathcal{M}_{2, 1} \times \mathcal{M}_{n - 3, 0} \subset \mathcal{M}_{n-1,1}$.
The $3$-dimensional Lorentzian metrics ${g_0}_{|E_{2, 1}}$ and $g_{|E_{2, 1}}$ 
are related by $g_{|E_{2, 1}} = {g_0}_{|E_{2, 1}} (B \cdot, \cdot)$, where $B$
is an automorphism of $E_{2, 1}$ of the form $k (\id + x)$, where $k$
is a positive real number and $x$ is an endomorphism of $E_{2, 1}$
satisfying $x^3=0$ but $x^2\neq 0$. 
Thus $g_0$ and $g$ are connected by a unique geodesic whose $E_{2, 1}$-part is of 
the form
\begin{equation} 
{g_t}_{|E_{2, 1}} = {g_0}_{|E_{2, 1}} (B_t\cdot, \cdot), \nonumber
\end{equation}
with $B_t = k^t \, {\rm exp} (t (x - \frac{1}{2} x^2)) = k^t \left(1 + t x +  
\frac{t (t - 1)}{2} x^2\right)$. 
\end{prop}

The rest of the paper is devoted to the proof of Proposition
\ref{proplorentz}. 

Recall that for  any $g$ and $g_0$ in $\mathcal{M}_{n - 1, 1}$, 
there exists
a uniquely defined
automorphism $A$ of $E$ --- with ${\rm det} A > 0$ --- such that $g = g_0 (A
\cdot, \cdot)$: $A = (\gamma^{-1})^* \gamma^{-1}$, for any
$\gamma\in\GL(E)$ such that $g = \gamma
\cdot g_0$ and $A$ is symmetric relative to both $g$ and $g_0$. 
Then $g_0$ can be joined with $g$ by a geodesic in $\mathcal{M}_{n-1,1}$ if and 
only if $A$ is of the form $A = {\rm exp} (a)$, for some
$g_0$-symmetric endomorphism $a$ of $E$, and the corresponding
geodesic is then the curve $g_t := g_0({\rm exp} (t a) \cdot, \cdot)$
for $t\in [0, 1]$. 

%In Case (i) of Proposition \ref{proplorentz}, $A$ is
%diagonalizable on an $(n - 2)$-dimensional subspace 
%$E_{n - 2}$ of
%$E$, entirely of space-type, where $A$ has only positive eigenvalues,
%and leaves invariant a $2$-dimensional complementary space  $E_{1, 1}$, which
%is orthogonal to $E_{n - 2, 0}$  and is of
%signature $(1, 1)$ for both $g_0$ and $g$. 
%Moreover, the restriction
%of $A$ to $E_{1, 1}$ can be of any kind described in Section~\ref{ssl2}. 
%In Case (ii) $A$ is diagonalizable on an
%$(n - 3)$-dimensional subspace $E_{n - 3}$ of
%$E$, entirely of space-type, where $A$ has only positive eigenvalues,
%and leaves invariant a $3$-dimensional complementary space  $E_{2, 1}$, which
%is orthogonal to $E_{n - 3, 0}$  and is of
%signature $(2, 1)$ for both $g_0$ and $g$. 
%Moreover, the restriction of $A$ to $E_{2, 1}$ is of the form $\id + x$,
%where $x$ satisfies $x^3 = 0$ but $x^2 \neq 0$ ($A_{|E_{2, 1}}$ has 
%then a $1$-dimensional, isotropic eigenspace). 

The proof of Proposition \ref{proplorentz} requires the spectral 
analysis of $A$.
For this purpose it is convenient to introduce a positive definite 
{\it Euclidean} inner product $(\cdot,\cdot)$ on $E$ such that 
$g_0 = (I \cdot, \cdot)$ where $I$ is of the form 
\begin{equation} \label{u} 
I = \id -  2 (u,\cdot) u, 
\end{equation}
for some element $u\in E$ such that $|u|^2 = 1$.
Here, and henceforth, $|\cdot|$ denotes the norm with respect to 
$(\cdot,\cdot)$.
For $g_0$ the vector $u$ is timelike with $g_0(u,u)=-1$. 
Conversely, any such $u$ determines a Euclidean inner product as above.
 
Then $g = g_0 (A \cdot, \cdot)$ can be written as $g  = (S \cdot,\cdot)$ for a 
uniquely defined $(\cdot,\cdot)$-symmetric automorphism $S$ of $E$
with exactly $n - 1$ positive and $1$ negative eigenvalues.

Conversely, for any such automorphism $S$, the inner product $g = (S \cdot,
\cdot)$ belongs to $\mathcal{M}_{n - 1, 1}$ with 
\begin{equation} 
A = I\hml S = I S. \nonumber
\end{equation}

The spectral decomposition of $S$ reads
\begin{equation} 
S = \lambda_0 \Pi_0 + \bigoplus_{r = 1}^{\ell}
  \lambda_j \Pi_r, \nonumber
\end{equation}
with $\lambda_0 < 0 < \lambda_1 <  \ldots \lambda_{\ell}$, where
$\Pi_j$ denotes the $(\cdot,\cdot)$-orthogonal projection onto the
$d_j$-dimensional eigenspace $E_j$ of $S$ corresponding to the 
eigenvalue $\lambda_j$, $j = 0, 1, \ldots, \ell$.
Note that $d_0 = 1$.

Via the decomposition $E = E_0 \oplus \bigoplus_{r = 1}^{\ell} E_r$ the
unit vector $u$ appearing in (\ref{u}) splits as
\begin{equation} u = u_0 + u_1 + \ldots + u_{\ell}. \nonumber\end{equation}

We denote by $\Delta$ the subset of $j\in\{ 0, 1, \ldots, \ell \}$
such that $u_j \neq 0$, and by $m$ the cardinality of $\Delta$. 
For each $j \in \Delta$ such that $d_j > 1$ we denote by $\tilde{E}_j$ the
$(\cdot,\cdot)$-orthogonal complement of $u_j$ in $E_j$.  
Let $\tilde{E}$ be the subspace
of $E$ defined by 
\begin{equation} \label{E0} 
\tilde{E}  := \bigoplus_{j \in \Delta,  d_j > 1} 
\tilde{E}_j \oplus \bigoplus_{j \notin \Delta} E_j,  
\end{equation}
and $W$ the $m$-dimensional subspace of $E$ defined by
\begin{equation} \label{W} 
W = \bigoplus_{j \in \Delta} \mathbb{R} \, u_j
\end{equation}
so that
\begin{equation} E = \tilde{E} \oplus W. \nonumber\end{equation}
Both $\tilde{E}$ and $W$ are left invariant by $A$, $I$, and $S$. 
The  sum is orthogonal  with respect to  $(\cdot,\cdot)$, $g_0$, and $g$. 

Note that if  $0 \notin \Delta$, i.~e.\ if $u_0 = 0$, then $\tilde{E}$ is of
signature $(n - m - 1, 1)$ and $W$ is of signature $(m,0)$, whereas,  if $0 \in
\Delta$, i.~e.\ if $u_0 \neq 0$, 
$W$ is of signature $(m - 1, 1)$ and $\tilde{E}$ is of signature $(n - m, 0)$ for $g$ (but $W$ is always of
signature $(m - 1, 1)$ for $g_0$, as $\tilde{E}$ is orthogonal to $u$).

Since $\tilde{E}$ is orthogonal to $u$, $I_{|\tilde{E}}  = \id$ and 
$A_{|\tilde{E}} = S_{|\tilde{E}}$. 
In particular, $A_{|\tilde{E}}$ is symmetric for $g_0$, $g$ and 
$(\cdot,\cdot)$ and its spectral decomposition coincides with the one of 
$S_{|\tilde{E}}$, given by (\ref{E0}), with eigenvalues $\lambda_j$ for each 
$j \notin \Delta$ and each $j \in \Delta$ with $d_j > 1$. 

The spectral study of $A$ is then reduced to the spectral study of $A
_{|W}$ and the latter is summarized by the following lemma.

\begin{lemma} \label{lemmaspec1} 
{\rm (i)} The characteristic polynomial $P$ of $A
_{|W}$ defined by $P (t) = {\rm det} (t\, \id - A_{|W})$ is given by
\begin{equation} \label{car1} 
P (t) = \prod_{j \in \Delta} (t -
  \lambda_j)  +  2 \sum_{j \in \Delta}
  \lambda_j |u_j|^2
  \prod_{k \in \Delta\setminus \{j\}} (t - \lambda_k). 
\end{equation}
In particular, the roots of $P$ are all distinct from the $\lambda
_j$, $j \in \Delta$. 

{\rm (ii)} For each real root $\mu$ of $P$ the corresponding eigenspace
 is the one-dimensional vector space generated by the element $v_{\mu}\in W$ 
defined by
\begin{equation} \label{vmu} v_{\mu}  =  \sum_{j \in \Delta} 
\frac{u_j}{\mu - \lambda_j}.  \end{equation}
Moreover, 
\begin{equation} \label{normv} g (v_{\mu}, v_{\mu}) =  \mu \, g_0 (v
_{\mu}, v_{\mu}) = - \frac{1}{2}
  \frac{P ' (\mu)}{Q (\mu)} \end{equation}
where $Q$ denotes the polynomial defined by $Q (t) = \prod_{j \in
  \Delta} (t - \lambda_j)$. In particular, $v_{\mu}$ is a null-vector
  --- for both $g$ and $g_0$ --- if and only if $\mu$ is a multiple
  root of $P$.
\end{lemma}
\begin{proof}
By definition, any $v \in W$ is of the form $v = \sum_{j \in
  \Delta} y_j u_j$, for real numbers $y_1, \ldots, y_m$, so that
\begin{equation} A v = I S v = \sum_{j \in \Delta} (\lambda_j y_j
 -  2 (Su, v) ) \, u_j. \nonumber\end{equation}
Note that $v$ is an eigenvector of $A_{|W}$ for some eigenvalue $\mu$ if 
and only if
\begin{equation} \label{propre} 
(\mu -  \lambda_j) \, y_j = - 2 (Su, v),
\end{equation}
for each $j \in \Delta$. 
It is easily checked  that $(Su, v)$ cannot be
equal 
to $0$ if $v \neq 0$. Indeed, suppose for a contradiction
that $v$ satisfies (\ref{propre}) with $(Su, v) = 0$ and $v \neq
0$. Since $v \neq 0$, one of the $y_j$, say $y_1$, is
nonzero, so that $\mu = \lambda_1$. This implies $\mu \neq  \lambda
_j$,  for $j \neq
1$, as the $\lambda_j$ are pairwise distinct. It follows that  $y_j = 0$ for
all $j \neq 1$, so that  $v = y_1 u_1$. Then  $(Su, v)
=  \lambda_1 y_1 |u_1|^2 \neq 0$ as  $y_1 \neq 0$, a contradiction.

In particular, this shows $\mu \neq \lambda_j$ for each $j\in\Delta$ 
so that we can write
\begin{equation} \label{v} 
v = - 2 (Su, v) \sum_{j \in \Delta} 
\frac{u_j}{\mu -  \lambda_j}. 
\end{equation}
Moreover, by computing $(Su, v) = (Sv, u)$ from (\ref{v}), we get
\begin{equation} \label{root} \sum_{j \in \Delta} \frac{\lambda_j |u_j|^2}{\mu
    -  \lambda_j} = - \frac{1}{2}.  \end{equation}
It follows that each eigenvalue of $A_{|W}$ is a root of the
polynomial $P$ defined by (\ref{car1}). Since $P$ is
monic and of degree $m$, it must coincide  with the
characteristic  polynomial of $A_{|W}$. 
We readily see from (\ref{car1}) that the roots of $P$ are distinct
from the $\lambda_j$ (recall that the latter are pairwise
distinct). From (\ref{v}) we immediately see that the eigenspace
corresponding to $\mu$ is generated by the vector $v_{\mu}$ defined
by (\ref{vmu}). 

Conversely, for each root $\mu$ of $P$ the vector $v
_{\mu}$ defined by (\ref{vmu}) is certainly an eigenvector of $A
_{|W}$ for the eigenvalue $\mu$. 

Since the roots of $P$ are distinct from the $\lambda_j$, $P$ can
also be expressed by
\begin{equation} \label{car2} 
\frac{P (t)}{Q (t)} = 1 +  2 \sum_{j \in \Delta}
  \frac{\lambda_j |u_j|^2}{t -  \lambda_j}, 
\end{equation}
where we put $Q (t) := \prod_{j \in \Delta} (t - \lambda_j)$.
Differentiating (\ref{car2}) at $t = \mu$, we get (\ref{normv}). 
It follows that $v_{\mu}$ is a null vector if and only if $P ' (\mu) = 0$, 
meaning that $\mu$ is a multiple root.
\end{proof}

For further use, we need more information about the sign of the
characteristic polynomial $P$ at $t = \lambda_j$, $j \in \Delta$, and
at $t = 0$. In the sequel, we use the notation $P (t_0) \equiv (-1)
^r$, for some integer $r$,  to mean that $P$ has the sign of $(-1)^r$ --- in particular is
not zero --- at $t = t_0$.

\begin{lemma} \label{lemmaspec2} {\rm (i)} If $0 \notin \Delta$, we
  re-label the $\lambda_j$ so that $\Delta = \{1, \ldots, m \}$, 
 and $0 < \lambda_1 < \ldots < \lambda_m$. We then have:
\begin{equation} \begin{split} & P (- \infty) \equiv P (\lambda_0)
 \equiv (-1)^m,\\ & P (0) \equiv (-1)^{m - 1},
 \\ & P (\lambda_j) \equiv (-1)^{m - j}, \ \ \ j = 1, \ldots, m.
 \end{split} \end{equation}
In particular, $P$ has then exactly $m$ distinct real roots $\mu_0 <
0 < \mu_1 < \ldots < \mu_{m - 1}$, with $\mu_0 \in (\lambda_0, 0)$
and $\mu_i \in (\lambda_i, \lambda_{i + 1})$, for $i = 1, \ldots, m -
1$. 

{\rm (ii)} If $0 \in \Delta$, we re-label the $\lambda_j$ so that
$\Delta = \{ 0, 1, \ldots, m - 1 \}$ and $\lambda_0 < 0 < \lambda_1
< \ldots < \lambda_{m - 1}$.We then have
\begin{equation} \begin{split} & P (- \infty) \equiv P (\lambda_0)
 \equiv P (0)  \equiv (-1)^m, \\ & P (\lambda
_j) \equiv (-1)^{m - j - 1}, \ \ \ j = 1, \ldots, m - 1.
 \end{split} \end{equation}
In particular, $P$ has then at least $(m - 2)$ distinct real roots $0
< \mu_1 < \ldots < \mu_{m - 2}$, with $\mu_i \in (\lambda_i,
\lambda_{i + 1})$, for $i = 1, \ldots, m-2$.
\end{lemma}

\begin{proof} 
Easy consequence of (\ref{car1}).
\end{proof}

We now consider the two cases when $0$ does or does not belong  to $\Delta$.

{\em Case 1}: $0 \notin \Delta$. 

According to Lemma~\ref{lemmaspec2} (i), $A_{|W}$ is diagonalizable 
(over $\mathbb{R}$) with one negative eigenvalue $\mu_0$ and $m - 1$ 
distinct positive eigenvalues. 
Moreover, we easily see from (\ref{normv}) that the $m$ corresponding
eigenvectors $v_{\mu}$, defined by (\ref{vmu}), are all spacelike.
On the other hand, $A_{|\tilde{E}}$ is also diagonalizable with one
negative eigenvalue, namely $\lambda_0$ --- whose 
eigenspace is $E_0$ ---  and $n - m - 1$ positive
eigenvalues. Denote by $E_{1, 1}$ the direct sum of $E_0$ and the
(one-dimensional) eigenspace of $\mu_0$, and by $E_{n - 2, 0}$ the
orthogonal complement of $E_{1, 1}$ for $g$ or $g_0$. 
Then, both $g$ and $g_0$ are of signature $(1, 1)$ on $E_{1, 1}$ and 
positive definite on $E_{n - 2, 0}$. 
Accordingly, $A$ splits as the sum of two operators $A
= A_{1, 1} \oplus A_{n - 2, 0}$, where $A_{1, 1}$ acts
trivially on $E_{n - 2, 0}$ and is diagonalizable, 
with negative eigenvalues on $E_{1,
  1}$, whereas $A_{n - 2, 0}$ acts trivially on $E_{1, 1}$ and
is positive definite, as well as $g_0$- and $g$-symmetric on $E_{n -  2,0} $. 
This can be interpreted as follows. Denote by
$\mathcal{M}_{1, 1}$ the space of Lorentzian inner products  of $E_{1, 1}$,
by $\mathcal{M}_{n - 2, 0}$ the space of positive definite inner
products of $E_{n - 2, 0}$. 
Then the product $\mathcal{M}_{1,1} \times \mathcal{M}_{n - 2, 0}$ is 
naturally embedded  as a totally geodesic submanifold of 
$\mathcal{M}_{n - 1, 1}$ and both 
$g = g_{|E_{1, 1}}\oplus g_{|E_{n-2,0}}$ and $g_0 = {g_0}
_{|E_{1, 1}}\oplus {g_0}_{|E_{n-2,0}}$  belong to it. 
In $\mathcal{M}_{n - 2, 0}$ any two elements, in particular $g_{|E
_{n - 2, 0}}$ and ${g_0}_{|E
_{n-2,0}}$,  are joined by a unique geodesic. The situation
  concerning $\mathcal{M}_{1, 1}$ has been explored in detail in the
  first part of this section. In the present case, $g_{|E_{1,
      1}}$ and ${g_0}_{|E_{1,
      1}}$ are related  by the automorphism $A_{|E_{1, 1}}$ which is
  diagonalizable with distinct negative eigenvalues, so that $g_{|E_{1,
      1}}$ and ${g_0}_{|E_{1,
      1}}$ cannot be linked by a geodesic. 

{\em Case 2}: $0 \in \Delta$.  

According to Lemma \ref{lemmaspec2} (ii), there exist at least $m-2$ 
distinct positive eigenvalues of $A_{|W}$, namely $0 < \mu_1 < \ldots < 
\mu_{m - 2}$. 
Then, either these eigenvalues are all simple roots of $P$, or one of them
--- and only one --- is a triple root.
The case that two of them are double roots is
impossible since, according to Lemma~\ref{lemmaspec1} (ii),  the
corresponding eigenvectors defined by
(\ref{vmu}) would then form an orthogonal pair of nonzero null
vectors in the Lorentzian space $(E, g)$. 

In the case when all $\mu_i$ are simple roots, we easily check by using
(\ref{normv}) that the corresponding eigenvectors are all spacelike. 
Denote by $E_{n - 2, 0}$
the direct sum of the corresponding eigenspaces and $\tilde{E}$, and by
$E_{1, 1} \subset W$ the orthogonal complement of $E_{n - 2, 0}$
 for $g$ or $g
_0$. Then, both $g$ and $g_0$ are positive definite on $E_{n - 2,
  0}$ and of signature $(1, 1)$ on $E_{1, 1}$. The situation is then
quite similar to the previous one, except that all cases considered in
Section \ref{ssl2} for $\mathcal{M}_{1, 1}$ may now happen, depending on 
whether the missing two roots of $P$ are complex conjugate, both positive 
(equal or distinct) or both negative (equal or distinct). 

It remains to consider the case that one of the $\mu_i$, say $\mu
_j := k > 0$, is a triple root of
$P$. Then, according to Lemma \ref{lemmaspec1} (iii), the
corresponding eigenvector $v_{\mu_j}$ is a null vector. 
Again, it is easily checked that the $v_{\mu_i}$, for $i \neq j$, are all
spacelike. 
Denote by $E_{n - 3, 0}$ the direct sum of the eigenspaces corresponding to 
the $\mu_i$, $i \neq j$, and $E^0$, and by $E_{2, 1} \subset W$ the 
orthogonal complement of $E_{n - 3, 0}$ for $g$ or $g_0$. 
Then, both $g$ and $g_0$ are positive definite on
$E_{n - 3, 0}$ and of signature $(2, 1)$ on $E_{2, 1}$. It follows
that $g$ and $g_0$ both belong to a same totally geodesic subspace
$\mathcal{M}_{2, 1} \times \mathcal{M}_{n - 3, 0}$. Moreover,
the restriction of $A$ to  $E_{2, 1}$, which relates $g|_{E_{2, 1}}$
and ${g_0}_{|E_{2, 1}}$,  is of the form $k (\id + x)$,
where $x$ is nilpotent and regular (this is because $\mu_j$ has no
other eigenvector than $v_{\mu_j}$). 
Now, $\id + x$ is the exponential
of $\id + x - \frac{x^2}{2}$, which is certainly symmetric for both $g
_0$ and $g$ (since $x = (\id + x) - \id$ is symmetric) and is the only
symmetric ``logarithm'' of $\id + x$. We thus get  a unique (null) geodesic
between ${g_0}_{|E_{2, 1}}$  and $g_{|E_{2, 1}}$ in $\mathcal{M}
_{2, 1}$, hence also between $g_0$ and $g$ in $\mathcal{M}_{n, 1}$. 

This completes  the proof of Proposition \ref{proplorentz}. 
$\hfill\Box$

%%%%%%%%%%%%%%%%%%%%%%%%%%%%%%%%%%%%%%%%%%%%%%%%%%%%%%%%%%%%%%%%%%%%%%%%%
%\bibliographystyle{amsplain}
%\bibliography{meine.bib,allg.bib}  
\providecommand{\bysame}{\leavevmode\hbox to3em{\hrulefill}\thinspace}

%%%%%%%%%%%%%%%%%%%%%%%%%%%%%%%%%%%%%%%%%%%%%%%%%%%%%%%%%%%%%%%%%%%%%%%%%

\end{document}

%% file: fig1c.tex
  \begin{pspicture}(-2,-2)(2,2)
    \psset{unit=0.66666666cm}
    %\psgrid[gridcolor=gray](-3,-3)(3,3)

    \psline[linestyle=dashed,dash=3pt 2pt](1.6,2.78)(1.6,-1.62)(-1.6,-2.78)
    \psline(1.6,1.25)(1.6,-1.25)

    \psclip{\pspolygon[linestyle=none](-3,-3)(-3,3)(3,3)(3,0)%
        (-0.6,0)(-0.6,-3)}
      \psellipse[linecolor=blue,linewidth=0.4pt,linestyle=dashed,%
        dash=3pt 2pt](0,0)(1.0,0.344)
    \endpsclip
    \psclip{\psframe[linestyle=none](-3,-2)(3,0)}
      \psellipse[linecolor=blue,linestyle=dashed,%
        dash=3pt 2pt](0,-2.2)(2.325,0.8)
    \endpsclip

    \parametricplot[linecolor=red,linestyle=dashed,dash=3pt 2pt]% 
    {-1.317}{1.317}{2.718 t exp 2.718 t neg exp add 2 div 0.4 sub
      2.718 t exp 2.718 t neg exp sub 2 div 0.2752 add 2.71 0.084 t
      mul 0.215 add exp mul} 

    \psellipse[linewidth=0.4pt,border=0.4pt](0,2.2)(2.09,0.72)
    \psellipse[linewidth=0.4pt,linestyle=dashed,dash=3pt 2pt]%
      (0,-2.2)(2.09,0.72)
    \psline[linewidth=0.4pt](-2,2)(-1.65,1.65)
    \psline[linewidth=0.4pt,linestyle=dashed,dash=3pt 2pt,%
      border=0.4pt](-1.65,1.65)(2,-2)
    \psline[linewidth=0.4pt,linestyle=dashed,dash=3pt 2pt,%
      border=0.4pt](-2,-2)(1.65,1.65)
    \psline[linewidth=0.4pt](1.65,1.65)(2,2)

    \psellipse[linecolor=blue,border=0.4pt](0,2.2)(2.325,0.8)
    \psclip{\psframe[linestyle=none](-0.6,-3)(3,0)}
      \psellipse[linecolor=blue,linewidth=0.4pt,border=0.4pt]%
        (0,0)(1.0,0.344)
    \endpsclip
    \psclip{\psframe[linestyle=none](-3,-3.5)(3,-1.9)}
      \psellipse[linecolor=blue](0,-2.2)(2.325,0.8)
    \endpsclip
    \parametricplot[linecolor=blue]{-1.44}{1.44}{2.718 t exp 2.718 t
      neg exp add 2 div 2.718 t exp 2.718 t neg exp sub 2 div}
    \psclip{\psframe[linestyle=none](-3,-3)(-1.6,3)}
      \parametricplot[linecolor=blue]{-1.44}{1.44}{2.718 t exp 2.718 t
        neg exp add 2 div neg 2.718 t exp 2.718 t neg exp sub 2 div}
    \endpsclip
    \psclip{\psframe[linestyle=none](-1.6,-3)(0,3)}
      \parametricplot[linecolor=blue,linestyle=dashed,dash=3pt 2pt]%
      {-1.44}{1.44}{2.718 t exp 2.718 t neg exp add 2 div neg 2.718 t
        exp 2.718 t neg exp sub 2 div} 
    \endpsclip

    \psline(1.6,2.78)(-1.6,1.62)(-1.6,-2.78)
    \psline[border=0.4pt](-1.6,1.5)(-1.6,-2.7)
    \psclip{\psframe[linestyle=none](1.5,-1.4)(1.7,1.4)}
      \parametricplot[linecolor=blue,border=0.4pt]{-1.44}{1.44}%
      {2.718 t exp 2.718 t neg exp add 2 div 2.718 t exp 2.718 t neg
        exp sub 2 div} 
    \endpsclip

    \parametricplot[linecolor=red]{-1.317}{1.317}{2.718 t exp 2.718 t
      neg exp add 2 div neg 0.4 add 2.718 t exp 2.718 t neg exp sub 2
      div 0.2752 sub 2.71 0.084 neg t mul 0.215 add exp mul}    
    \psclip{%
        \pscustom[linestyle=none]{%
          \psframe[linestyle=none](-1.5,1.5)(-1.0,0.8)
          \psframe(-1.55,-2.75)(-0.4,-0.4)}}
      \parametricplot[linecolor=red,border=0.5pt]{-1.317}{1.317}{2.718 t exp 2.718 t
        neg exp add 2 div neg 0.4 add 2.718 t exp 2.718 t neg exp sub 2
        div 0.2752 sub 2.71 0.084 neg t mul 0.215 add exp mul}    
    \endpsclip

    \rput[l](1.7,0){$E$}
    \rput[br](-2.4,-2.2){\blue $\MM^{1,1}_1$}
    \psset{unit=1cm}
  \end{pspicture}

%% file: fig1a.tex
  \begin{pspicture}(-2,-2)(2,2)
    \psset{unit=0.66666666cm}
    %\psgrid[gridcolor=gray](-3,-3)(3,3)

    \psclip{\pspolygon[linestyle=none](-3,-3)(-3,1)(-1.1,1)%
        (-1.1,0)(1.4,0)(1.4,2)(3,2)(3,-3)}
      \psline(-1.6,0.208)(2.0,1.168)
    \endpsclip
    \psclip{\pspolygon[linestyle=none](-1.1,1)(-1.1,0)(1.4,0)(1.4,2)}
      \psline[linestyle=dashed,dash=3pt 2pt](-1.6,0.208)(2.0,1.168)
    \endpsclip

    \psclip{\pspolygon[linestyle=none](-1.8,-0.48)(1.8,0.48)%
        (2.0,2.168)(-1.6,0.208)}
      \rput{15}(0,0){\psellipse[linecolor=red,linestyle=dashed,%
        dash=3pt 2pt](1.1,0.344)}
    \endpsclip

    \psclip{\psframe[linestyle=none](-3,3)(3,0)}
      \psellipse[linecolor=blue,linewidth=0.4pt,linestyle=dashed,%
        dash=3pt 2pt](0,0)(1.0,0.344)
    \endpsclip
    \psclip{\psframe[linestyle=none](-3,-2)(3,0)}
      \psellipse[linecolor=blue,linestyle=dashed,%
        dash=3pt 2pt](0,-2.2)(2.325,0.8)
    \endpsclip

    \psellipse[linewidth=0.4pt,border=0.4pt](0,2.2)(2.09,0.72)
    \psellipse[linewidth=0.4pt,linestyle=dashed,dash=3pt 2pt]%
      (0,-2.2)(2.09,0.72)
    \psline[linewidth=0.4pt](-2,2)(-1.65,1.65)
    \psline[linewidth=0.4pt,linestyle=dashed,dash=3pt 2pt,%
      border=0.4pt](-1.65,1.65)(2,-2)
    \psline[linewidth=0.4pt,linestyle=dashed,dash=3pt 2pt,%
      border=0.4pt](-2,-2)(1.65,1.65)
    \psline[linewidth=0.4pt](1.65,1.65)(2,2)

    \psellipse[linecolor=blue,border=0.4pt](0,2.2)(2.325,0.8)
    \psclip{\psframe[linestyle=none](-3,-3)(3,0)}
      \psellipse[linecolor=blue,linewidth=0.4pt,border=0.4pt]%
        (0,0)(1.0,0.344)
    \endpsclip
    \psclip{\psframe[linestyle=none](-3,-3.5)(3,-1.9)}
      \psellipse[linecolor=blue](0,-2.2)(2.325,0.8)
    \endpsclip
    \psclip{\pspolygon[linestyle=none](-2.0,-1.168)(1.6,-0.208)%
        (1.8,0.4)(-1.8,-0.4)}
      \parametricplot[linecolor=blue,linestyle=dashed,dash=3pt 2pt]%
      {-1.44}{1.44}{2.718 t exp 2.718 t neg exp add 2 div 2.718 t exp
        2.718 t neg exp sub 2 div} 
      \parametricplot[linecolor=blue,linestyle=dashed,dash=3pt 2pt]%
      {-1.44}{1.44}{2.718 t exp 2.718 t neg exp add 2 div neg 2.718 t
        exp 2.718 t neg exp sub 2 div} 
    \endpsclip
    \psclip{%
      \pscustom[linestyle=none]{%
        \psframe(-3,-3)(3,3)
        \pspolygon(-2.0,-1.168)(1.6,-0.208)(1.8,0.4)(-1.8,-0.4)}}
      \parametricplot[linecolor=blue]{-1.44}{1.44}{2.718 t exp 2.718 t
        neg exp add 2 div 2.718 t exp 2.718 t neg exp sub 2 div}
      \parametricplot[linecolor=blue]%
      {-1.44}{1.44}{2.718 t exp 2.718 t neg exp add 2 div neg 2.718 t
        exp 2.718 t neg exp sub 2 div} 
    \endpsclip
    \psline[border=0.4pt](-2.0,-1.168)(1.6,-0.208)
    \psline(-1.6,0.208)(-2.0,-1.168)
    \psline(1.6,-0.208)(2.0,1.168)
    \psclip{\psframe[linestyle=none](1.6,1.2)(1.2,0.8)}
      \parametricplot[linecolor=blue,border=0.4pt]{-1.44}{1.44}%
      {2.718 t exp 2.718 t neg exp add 2 div 2.718 t exp 2.718 t neg
        exp sub 2 div} 
    \endpsclip
    \psclip{\psframe[linestyle=none](-1.2,0.2)(-0.8,0.6)}
      \parametricplot[linecolor=blue,border=0.4pt]{-1.44}{1.44}%
      {2.718 t exp 2.718 t neg exp add 2 div neg 2.718 t exp 2.718 t
      neg exp sub 2 div} 
    \endpsclip

    \psclip{\pspolygon[linestyle=none](-1.8,-0.48)(1.8,0.48)%
        (1.6,-0.208)(-2.0,-1.168)}
      \rput{15}(0,0){\psellipse[linecolor=red](1.1,0.344)}
    \endpsclip
    
    \rput[tl](1.7,-0.2){$E$}
    \rput[br](-2.4,-2.2){\blue $\MM^{1,1}_1$}
    \psset{unit=1cm}
  \end{pspicture}

%% file: fig1b.tex
  \begin{pspicture}(-2,-2)(2,2)
    \psset{unit=0.66666666cm}
    %\psgrid[gridcolor=gray](-3,-3)(3,3)

    \psclip{\pspolygon[linestyle=none](-3,-2.4)(-2.1,-2.4)(2.1,2.4)%
        (2.1,3)(-3,3)}
      \pspolygon[linestyle=dashed,dash=3pt 2pt]%
        (2.581,3)(1.419,1)(-2.581,-3)(-1.419,-1)
    \endpsclip

    \psclip{\psframe[linestyle=none](-3,3)(3,0)}
      \psellipse[linecolor=blue,linewidth=0.4pt,linestyle=dashed,%
        dash=3pt 2pt](0,0)(1.0,0.344)
    \endpsclip
    \psclip{\psframe[linestyle=none](-3,-2)(3,0)}
      \psellipse[linecolor=blue,linestyle=dashed,%
        dash=3pt 2pt](0,-2.2)(2.325,0.8)
    \endpsclip

    \psline[linecolor=red,linestyle=dashed,dash=3pt 2pt]%
      (2.25,2.41)(-1.82,-1.66)
    \psline[linecolor=red](2.25,2.41)(2.05,2.21)  

    \psellipse[linewidth=0.4pt,border=0.4pt](0,2.2)(2.09,0.72)
    \psellipse[linewidth=0.4pt,linestyle=dashed,dash=3pt 2pt]%
      (0,-2.2)(2.09,0.72)
    \psline[linewidth=0.4pt](-2,2)(-1.65,1.65)
    \psline[linewidth=0.4pt,linestyle=dashed,dash=3pt 2pt,%
      border=0.4pt](-1.65,1.65)(2,-2)
    \psline[linewidth=0.4pt,linestyle=dashed,dash=3pt 2pt,%
      border=0.4pt](-2,-2)(1.65,1.65)
    \psline[linewidth=0.4pt](1.65,1.65)(2,2)

    \psellipse[linecolor=blue,border=0.4pt](0,2.2)(2.325,0.8)
    \psclip{\psframe[linestyle=none](-3,-3)(3,0)
        \pspolygon[linestyle=none](1.419,1)(-2.581,-3)(-2.25,-2.41)(1.82,1.66)}
      \psellipse[linecolor=blue,linestyle=dashed,dash=1pt 2pt,%
        linewidth=0.4pt,border=0.4pt](0,0)(1.0,0.344)
    \endpsclip
    \psclip{\psframe[linestyle=none](-3,-3.5)(3,-1.9)
        \pspolygon[linestyle=none](1.419,1)(-2.581,-3)(-2.25,-2.41)(1.82,1.66)}
      \psellipse[linecolor=blue,linestyle=dashed,dash=3pt 2pt]%
        (0,-2.2)(2.325,0.8)
    \endpsclip
    \psclip{\pspolygon[linestyle=none](-3,-3)(-3,0)(-0.2,0)(-0.2,-1)%
        (0.1,-1)(0.1,0)(3,0)(3,-3)}
      \psellipse[linecolor=blue,linewidth=0.4pt,border=0.4pt]%
        (0,0)(1.0,0.344)
    \endpsclip
    \psclip{\pspolygon[linestyle=none](-3,-3.5)(-3,-1.9)(-2.15,-1.9)%
        (-2.25,-2.7)(-2.1,-2.7)(-2.1,-1.9)(3,-1.9)(3,-3.5)}
      \psellipse[linecolor=blue](0,-2.2)(2.325,0.8)
    \endpsclip
    \parametricplot[linecolor=blue]{-1.44}{1.44}{2.718 t exp 2.718 t
      neg exp add 2 div 2.718 t exp 2.718 t neg exp sub 2 div}
    \parametricplot[linecolor=blue]{-1.44}{1.44}{2.718 t exp 2.718 t
      neg exp add 2 div neg 2.718 t exp 2.718 t neg exp sub 2 div}

    \psline[linecolor=red](-2.25,-2.41)(1.82,1.66)
    \psline[linecolor=red,border=0.4pt](-2.20,-2.36)(-0.84,-1.0)
    \psline[linecolor=red,border=0.4pt](-0.04,-0.20)( 0.76, 0.6)

    \psclip{\pspolygon[linestyle=none](-3,-4)(-3,-2.4)(-2.1,-2.4)%
        (2.1,2.4)(2.1,4)(3,4)(3,-4)}
      \pspolygon(2.581,3)(1.419,1)(-2.581,-3)(-1.419,-1)
    \endpsclip
    \psclip{\pspolygon[linestyle=none](-2.2,-3)(-2.2,1.6)(3,1.6)(3,-3)}
      \psline[border=0.4pt](2.581,3)(1.419,1)(-2.581,-3)
    \endpsclip
    \psclip{\psframe[linestyle=none](1.8,2.45)(2.15,2.8)}
      \psellipse[linecolor=blue,border=0.4pt](0,2.2)(2.325,0.8)
    \endpsclip
    \psclip{\psframe[linestyle=none](-2,-1.6)(-1,-0.4)}
      \parametricplot[linecolor=blue,border=0.4pt]{-1.44}{1.44}%
      {2.718 t exp 2.718 t neg exp add 2 div neg 2.718 t exp 2.718 t
        neg exp sub 2 div} 
    \endpsclip

    \rput[tl](2.6,2.8){$E$}
    \rput[br](-2.4,-2.2){\blue $\MM^{1,1}_1$}
    \psset{unit=1cm}
  \end{pspicture}

%% file: fig2.tex
  \begin{pspicture}(-2,-2)(2,2)
    \psset{unit=0.66666666cm}
    %\psgrid[gridcolor=gray](-3,-3)(3,3)

    \psclip{\psframe[linestyle=none](-3,3)(3,0)}
      \psellipse[linecolor=blue,linewidth=0.4pt,linestyle=dashed,%
        dash=3pt 2pt](0,0)(1.0,0.344)
    \endpsclip
    \psclip{\psframe[linestyle=none](-3,-2)(3,0)}
      \psellipse[linecolor=blue,linestyle=dashed,%
        dash=3pt 2pt](0,-2.2)(2.325,0.8)
    \endpsclip

    \psellipse[linecolor=blue,border=0.4pt](0,2.2)(2.325,0.8)
    \psclip{\psframe[linestyle=none](-3,-3)(3,0)
        \pspolygon[linestyle=none](1.419,1)(-2.581,-3)(-2.25,-2.41)(1.82,1.66)}
      \psellipse[linecolor=blue,linestyle=dashed,dash=1pt 2pt,%
        linewidth=0.4pt,border=0.4pt](0,0)(1.0,0.344)
    \endpsclip
    \psclip{\psframe[linestyle=none](-3,-3.5)(3,-1.9)
        \pspolygon[linestyle=none](1.419,1)(-2.581,-3)(-2.25,-2.41)(1.82,1.66)}
      \psellipse[linecolor=blue,linestyle=dashed,dash=3pt 2pt]%
        (0,-2.2)(2.325,0.8)
    \endpsclip
    \psclip{\pspolygon[linestyle=none](-3,-3)(-3,0)(-0.2,0)(-0.2,-1)%
        (0.1,-1)(0.1,0)(3,0)(3,-3)}
      \psellipse[linecolor=blue,linewidth=0.4pt,border=0.4pt]%
        (0,0)(1.0,0.344)
    \endpsclip
    \psclip{\pspolygon[linestyle=none](-3,-3.5)(-3,-1.9)(-2.15,-1.9)%
        (-2.25,-2.7)(-2.1,-2.7)(-2.1,-1.9)(3,-1.9)(3,-3.5)}
      \psellipse[linecolor=blue](0,-2.2)(2.325,0.8)
    \endpsclip
    \parametricplot[linecolor=blue]{-1.44}{1.44}{2.718 t exp 2.718 t
      neg exp add 2 div 2.718 t exp 2.718 t neg exp sub 2 div}
    \parametricplot[linecolor=blue]{-1.44}{1.44}{2.718 t exp 2.718 t
      neg exp add 2 div neg 2.718 t exp 2.718 t neg exp sub 2 div}

    \psline[linecolor=red](-2.25,-2.41)(1.82,1.66)
    \psline[linecolor=red,border=0.4pt](-2.20,-2.36)(-0.84,-1.0)
    \psline[linecolor=red,border=0.4pt](-0.04,-0.20)( 0.76, 0.6)
    \psline[linecolor=red](2.0,-2.51)(-1.42,1.01)
    
    \psclip{\psframe[linestyle=none](1.8,2.45)(2.15,2.8)}
      \psellipse[linecolor=blue,border=0.4pt](0,2.2)(2.325,0.8)
    \endpsclip
    \psclip{\psframe[linestyle=none](-2,-1.6)(-1,-0.4)}
      \parametricplot[linecolor=blue,border=0.4pt]{-1.44}{1.44}%
      {2.718 t exp 2.718 t neg exp add 2 div neg 2.718 t exp 2.718 t
        neg exp sub 2 div} 
    \endpsclip

    \psecurve[linewidth=0.4pt]{->}(3.5,-1)(2.8,0.5)(-.3,0.9)(-1,0)
    \psecurve[linewidth=0.4pt]{->}(2.5,2)(2.8,-0.5)(0,-2)(-1,-1)

    \psdots(-.173,-.325)
    \rput(-0.3,-0.8){$-I$}
    \psdots(.173,.325)
    \rput(0.4,0.7){$I$}
    \rput[br](-2.4,2.2){\blue $\MM^{1,1}_1$}
    \rput(2.8,0.25){unreachable} 
    \rput(2.8,-0.25){points}
    \psset{unit=1cm}
  \end{pspicture}